\documentclass[11pt,reqno,twoside]{amsart}

\usepackage{amsmath, amssymb, amscd, amsthm}

\theoremstyle{plain}
\newtheorem{thm}{Theorem}[section]

\newtheorem{lem}[thm]{Lemma}
\newtheorem{prop}[thm]{Proposition}
\newtheorem{Def}[thm]{Definition}

\theoremstyle{remark}
\newtheorem{rem}[thm]{Remark}
\newcounter{sspar}[subsection]
\renewcommand\thesspar{(\thesubsection.\arabic{sspar})}
    {\par\ \newline
     \vskip-\baselineskip\vskip.1truecm
     \noindent\refstepcounter{sspar}
     \noindent\textbf{\thesspar} \ignorespaces}
    {\vskip-\baselineskip
    \ignorespaces}
    {\refstepcounter{sspar}
     \textup{\textbf{\thesspar}} \ignorespaces}
    {\vskip-\baselineskip
    \ignorespaces}
\newcommand{\R}{\mathbb R}
\newcommand{\N}{\mathbb N}
\newcommand{\C}{\mathbb C}
\newcommand{\Z}{\mathbb Z}
\newcommand{\al}{\alpha}
\newcommand{\be}{\beta}

\newcommand{\Ga}{\Gamma}
\newcommand{\de}{\delta}
\newcommand{\De}{\Delta}
\newcommand{\eps}{\varepsilon}
\newcommand{\si}{\sigma}

\newcommand{\Te}{\Theta}
\newcommand{\la}{\lambda}
\newcommand{\La}{\Lambda}
\newcommand{\ph}{\varphi}
\newcommand{\Ups}{\Upsilon}

\newcommand{\Om}{\Omega}
\newcommand{\thet}{\vartheta}
\newcommand{\lt}{\ell^2}
\newcommand{\inprod}[2]{\langle #1,#2 \rangle}
\newcommand{\tensor}{\otimes}
\newcommand{\dirint}{\sideset{}{^\oplus}\int\limits_{0}^{\infty}}

\newcommand{\F}[5]{\,_{#1}F_{#2} \left( \genfrac{.}{.}{0pt}{}{#3}{#4}
\ ;#5 \right)}
\newcommand{\hf}{\frac{1}{2}}
\newcommand{\su}{\mathfrak{su}}

\newcommand{\calS}{\mathcal{S}}
\newcommand{\Res}[1]{\underset{#1}{\mathrm{Res}}}
\newcommand{\vect}[2]{ \begin{pmatrix} #1 \\ #2 \end{pmatrix} }

\newcommand{\D}{\mathcal D}
\numberwithin{equation}{section}

\begin{document}

\date{February 20, 2003}
\title[Continuous Hahn functions]{Continuous Hahn functions as Clebsch-Gordan coefficients}
\author{Wolter Groenevelt}
\author{Erik Koelink}
\address{
Technische Universiteit Delft, ITS-TWA \\
Postbus 5031, 2600 GA Delft, The Netherlands}
\email{W.G.M.Groenevelt@its.tudelft.nl  \\ H.T.Koelink@its.tudelft.nl}
\author{Hjalmar Rosengren}
\address{
Department of Mathematics, Chalmers University of Technology and G\"oteborg University\\
SE-412 96 G\"oteborg, Sweden}
\email{hjalmar@math.chalmers.se}

\dedicatory{Dedicated to Mizan Rahman}

\begin{abstract}
An explicit bilinear generating function for Meixner-Pollaczek polynomials is proved. This formula involves continuous dual Hahn polynomials, Meixner-Pollaczek functions, and non-polynomial $_3F_2$-hypergeometric functions that we consider as continuous Hahn functions. An integral transform pair with continuous Hahn functions as kernels is also proved. These results have an interpretation for the tensor product decomposition of a positive and a negative discrete series representation of $\su(1,1)$ with respect to hyperbolic bases, where the Clebsch-Gordan coefficients are continuous Hahn functions.
\end{abstract}
\maketitle

\section{Introduction}
The results and techniques in this paper are mainly analytic in nature, but they are motivated by a Lie algebraic problem. As is well known, many polynomials in the Askey-scheme of orthogonal polynomials of hypergeometric type, see \cite{KS}, have an interpretation in the representation theory of Lie groups and Lie algebras, see e.g.~ Vilenkin and Klimyk \cite{VK} and Koornwinder \cite{Koo2}. The Askey-scheme can be extended to families of unitary integral transforms with a hypergeometric kernel. Many of these kernels also admit group theoretic interpretations. For example the Jacobi functions, which can be considered as a non-polynomial extension of the Jacobi polynomials and are given explicitly by a certain $_2F_1$-hypergeometric function, have an interpretation as matrix elements for irreducible representations of the Lie group $SU(1,1)$. The Jacobi function is the kernel in the Jacobi integral transform, which can be found by spectral analysis of the hypergeometric differential operator. For an overview of Jacobi functions in representation theory, we refer to the survey paper \cite{Koo} by Koornwinder. 

In this paper we give a generalization of the Jacobi functions. We consider the tensor product of a positive and a negative discrete series representation of the Lie algebra $\su(1,1)$. The Clebsch-Gordan coefficients for the hyperbolic basisvectors turn out to be a certain type of non-polynomial $_3F_2$-hypergeometric functions, which we call continuous Hahn functions. We show that the continuous Hahn functions are the kernel in an integral transform, that generalizes the Jacobi function transform. We emphasize that the main part (sections \ref{sec2} and \ref{sec:Hahn}) of this paper is analytic in nature, and that the Lie algebraic interpretation is mainly restricted to section \ref{sec:su11}.

The Lie algebra $\su(1,1)$ is generated by the three elements $H$, $B$ and $C$. There are four classes of irreducible unitary representations for $\su(1,1)$: discrete series, i.e.~ the positive and the negative discrete series representations, and continuous series, i.e.~ the principal unitary series and the complementary series representations. There are three kinds of basis elements on which the various representations can act: the elliptic, the parabolic and the hyperbolic basis elements. These three elements are related to conjugacy classes of the group $SU(1,1)$. We consider the tensor product of a positive and a negative discrete series representation, which decomposes into a direct integral over the principal unitary series representations. Under certain condition discrete terms can appear. The Clebsch-Gordan coefficients for the standard (elliptic) basis vectors are continuous dual Hahn polynomials.
We compute the Clebsch-Gordan coefficients for the hyperbolic basis vectors, which are non-polynomial extensions of the continuous (dual) Hahn polynomials, and are therefore called continuous Hahn functions. For the Clebsch-Gordan coefficients for the elliptic and parabolic basis, we refer to \cite{GK}, respectively \cite{BW}, \cite{WG}.

The explicit expressions for the Clebsch-Gordan coefficients as $_3F_2$-series are not new, they are found by Mukunda and Radhakrishnan in \cite{MuR}. However not much seems to be known about the generalized orthogonality properties of the continuous Hahn functions, i.e.~ they form the kernel in a unitary integral tranform (the continuous Hahn transform). Using the Lie algebraic interpretation of the continuous Hahn functions, we can compute formally the inverse of the continuous Hahn integral transform. In section \ref{sec:Hahn} we give an analytic proof for the integral transform pair. 

The method we use to compute the Clebsch-Gordan coefficients is based on an idea by Granovskii and Zhedanov \cite{GZ}. The idea is to consider a  self-adjoint Lie algebra element $X_a= -aH+B-C$, $a\in \R$. The action of $X_a$ in an irreducible representation gives a difference equation, for which the (generalized) eigenvectors can be expressed in terms of special functions and the standard basis vectors. The Clebsch-Gordan coefficients for the eigenvectors can be calculated using properties of the special functions. In \cite{VdJ} and \cite{KJ1} Van der Jeugt and the second author considered the action of $X_a$ in tensor products of positive discrete series representations of $\su(1,1)$ to find convolution formulas for orthogonal polynomials. In \cite{GK} the action of $X_a$ in the tensor product of a positive and a negative discrete series representation is investigated for $|a|>1$ (the elliptic case). This leads to a bilinear summation formula for Meixner polynomials \cite[Thm.3.6]{GK}. In this paper we consider the case $|a|<1$ (the hyperbolic case). 

The plan of the paper is as follows. In section \ref{sec:functions} we introduce the special functions we need in this paper, and give some properties of these functions.

In section \ref{sec2} we prove a bilinear summation formula for Meixner-Pollaczek polynomials by series manipulations. As a result we find a certain type of $_3F_2$-functions, which are the continuous Hahn functions. The summation formula is used in section \ref{ssec:H} to compute the Clebsch-Gordan coefficients for the hyperbolic bases.

In section \ref{sec:su11} we consider the tensor product of a positive and a negative discrete series representation of the Lie algebra $\su(1,1)$. First  we recall the basic properties of $\su(1,1)$ and its irreducible unitary representations in section \ref{ssec:su11}. Then in section \ref{ssec:H} we diagonalize $X_a$, $|a|<1$, in the various irreducible representations. This leads to generalized eigenvectors of $X_a$, which can be considered as hyperbolic basis vectors.

For the discrete series representations, the overlap coefficients for the eigenvectors and the standard (elliptic) basisvectors are Meixner-Pollaczek polynomials, cf.~ \cite[\S7]{Koo2}. For the continuous series, the overlap coefficients are Meixner-Pollaczek functions. This follows from the spectral analysis of a doubly infinite Jacobi operator, which is carried out by Masson and Repka \cite[\S3.3]{MR} and Koelink \cite[\S4.4.11]{Koe}. It turns out that the spectral projection of the Jacobi operator is on a $2$-dimensional space of generalized eigenvectors. So the eigenvectors of $X_a$, $|a|<1$, in the continuous series representations are $2$-dimensional, and we find two linearly independent Meixner-Pollaczek functions as overlap coefficients. To determine the Clebsch-Gordan coefficients for the hyperbolic bases, we use the bilinear summation formula from section \ref{sec2}. This leads to a pair of continuous Hahn functions as Clebsch-Gordan coefficients. By formal calculations we find an integral transform pair, with a pair of continuous Hahn functions as a kernel. To give a rigorous proof of the integral transform pair, we show that the continuous Hahn functions are eigenfunctions of a difference operator $\La$. To find this operator $\La$ we realize $H$, $B$ and $C$ as difference operators acting on polynomials, using the difference equation for the Meixner-Pollaczek polynomials. Then $\La$ is a restriction of the Casimir operator in the tensor product.

The spectral analyis of this difference operator is carried out in section \ref{sec:Hahn}. A main problem with spectral analysis of a difference operator is finding the right eigenfunctions. This is because an eigenfunction multiplied by a periodic function is again an eigenfunction. Our choice of the periodic function is mainly motivated by the Lie algebraic interpretation of the eigenfunctions. Using asymptotic methods, we find a spectral measure for the difference operator. This leads to an integral transform with a pair of continuous Hahn functions as a kernel. We call this the continuous Hahn integral transform. \\

\emph{Notations.}
We denote for a function $f: \C \rightarrow  \C$
\[
f^*(x) = \overline{f(\overline{x})}.
\]
If $d\mu(x)$ is a positive measure, we use the notation $d\mu^\hf(x)$ for the positive measure with the property 
\[
d(\mu^\hf \times \mu^\hf)(x,x) = d\mu(x).
\]
The hypergeometric series is defined by
\[
\F{p}{q}{a_1, \ldots, a_{p}} {b_1, \ldots, b_q}{z} = \sum_{n=0}^\infty \frac{ (a_1)_n \ldots (a_{p})_n }{ (b_1)_n \ldots (b_q)_n } \frac{z^n}{n!},
\]
where $(a)_n$ denotes the Pochhammer symbol, defined by
\[
(a)_n = \frac{\Ga(a+n)}{\Ga(a)}=a(a+1)(a+2) \ldots (a+n-1), \qquad n \in \Z_{\geq 0}.
\]

\emph{Acknowledgements.} We thank Ben de Pagter for useful discussions.\\

\emph{Dedication.} We gladly dedicate this paper to Mizan Rahman who, with his
unsurpassed mastery in dealing with ($q$-)series and his insight
in the structures of formulas, has pushed the subject of ($q$)-special functions much further. We are also grateful 
to Mizan Rahman for his interest in our work, and for 
his willingness to help others in solving problems in this field.

\section{Orthogonal polynomials and functions} \label{sec:functions}
In this section we recall some properties of the orthogonal polynomials and functions which we need in this paper.\\

\textbf{Continuous dual Hahn polynomials.}
The Wilson polynomials, see Wilson \cite{Wil} or \cite[\S3.8]{AAR}, are $_4F_3$-hypergeometric polynomials on top of the Askey-scheme of hypergeometric polynomials, see Koekoek and Swarttouw \cite{KS}. The continuous dual Hahn polynomials are a three-parameter subclass of the Wilson polynomials, and are defined by 
\begin{equation} \label{def:cont dHahn}
s_n(y;a,b,c) = (a+b)_n (a+c)_n\F{3}{2}{-n,a+ix,a-ix}{a+b,a+c}{1},
\quad x^2=y, \quad (n \in \Z_{\geq 0}).
\end{equation}
For real parameters $a$, $b$, $c$, with $a+b$, $a+c$, $b+c$ positive,
the continuous dual Hahn polynomials are orthogonal with respect to a
positive measure, supported on a subset of $\R$. The orthonormal
continuous dual Hahn polynomials are defined by
\[
S_n(y;a,b,c) = \frac{(-1)^n s_n(y;a,b,c)} {\sqrt{n!(a+b)_n
(a+c)_n (b+c)_n}}\ .
\]
By Kummer's transformation, see e.g. \cite[Cor. 3.3.5]{AAR}, the
polynomials $s_n$ and $S_n$ are symmetric in $a$, $b$ and $c$.
Without loss of generality we assume that $a$ is the smallest of the
real parameters $a$, $b$ and $c$. Let $d\mu(\cdot;a,b,c)$ be the
measure defined by
\begin{align*}
\int_{\R} f(y) d&\mu(y;a,b,c)  =
 \frac{1}{2\pi}\int_0^\infty \left|
\frac{\Gamma(a+ix)\Gamma(b+ix)\Gamma(c+ix)}{\Gamma(2ix)} \right|^2
\frac{f(x^2)}{\Gamma(a+b) \Gamma(a+c) \Gamma(b+c)}\,d x\\ + &
\frac{\Gamma(b-a) \Gamma(c-a)}{\Gamma(-2a) \Gamma(b+c)} \sum_{k=0}^K
(-1)^k \frac{(2a)_k (a+1)_k (a+b)_k (a+c)_k} {(a)_k (a-b+1)_k (a-c+1)_k
k!} f(-(a+k)^2),
\end{align*}
where $K$ is the largest non-negative integer such that $a+K<0$. In
particular, the measure $d\mu(\cdot;a,b,c)$ is absolutely continuous
if $a\geq0$. The measure is positive under the conditions $a+b>0$,
$a+c>0$ and $b+c>0$. Then the polynomials $S_n(y;a,b,c)$ are
orthonormal with respect to the measure $d\mu(y;a,b,c)$.\\

\textbf{Meixner-Pollaczek polynomials.} The Meixner-Pollaczek polynomials, see \cite{KS}, \cite[Ex.6.37]{AAR}, are a two-parameter subclass of the Wilson polynomials, and are defined by
\begin{equation} \label{def:MP pol}
p_n^{(\lambda)}(x;\phi) = \frac{(2\lambda)_n}{n!} e^{in \phi}
\F{2}{1}{-n,\lambda+ix}{2\lambda}{1-e^{-2i\phi}}, \quad (n \in \Z_{\geq 0}).
\end{equation}
For $\lambda>0$ and $0<\phi<\pi$, these are orthogonal polynomials with respect to a positive measure on $\mathbb{R}$. The orthonormal Meixner-Pollaczek polynomials
\[
P_n(x) = P_n^{(\lambda)}(x;\phi) = \sqrt{\frac{n!} {(2
\lambda)_n}} p_n^{(\lambda)}(x;\phi)
\]
satisfy the following three-term recurrence relation
\begin{equation} \label{MP rec}
2x \sin{\phi}\, P_n(x) = \alpha_n P_{n+1}(x) - 2(n+\lambda)
\cos{\phi}\,P_n(x) + \alpha_{n-1} P_{n-1}(x),
\end{equation}
where
\[
\alpha_n = \sqrt{(n+1)(n+2 \lambda)}.
\]
The Meixner-Pollaczek polynomials also satisfy the difference equation
\begin{equation} \label{diff.eq.MPpol}
\begin{split}
e^{i\phi}(\la-ix) y(x+i) + 2i[x\cos\phi - (n+\la)\sin\phi] y(x)& - e^{-i\phi}(\la+ix) y(x-i) =0,\\
y(x) &= P_n^{(\lambda)}(x;\phi).
\end{split}
\end{equation}
Define
\[
w^{(\lambda)}(x;\phi) = \frac{1} {2 \pi} (2 \sin{\phi})^{2\lambda}e^{(2
\phi - \pi)x}\frac{ |\Gamma(\lambda+ix)|^2}{\Ga(2\la)},
\]
then the orthonormality relation reads
\[
\int_{\R} P_{m}^{(\lambda)}(x;\phi) P_{n}^{(\lambda)}(x;\phi)
w^{(\lambda)}(x;\phi) d x = \delta_{mn}.
\]

\textbf{Meixner-Pollaczek functions.}
The Meixner-Pollaczek functions $u_n$, $n \in \Z$, can be considered
as non-polynomial extensions of the Meixner-Pollaczek polynomials,
see Masson and Repka \cite{MR} and Koelink \cite[\S4.4]{Koe}. The Meixner-Pollaczek functions are defined by
\begin{equation} \label{def:MPfunctions}
\begin{split}
u_n(x;\la,\eps,\phi) =  (2i\sin \phi)^{-n} &\frac{
\sqrt{ \Ga(n+1+\eps+\la) \Ga(n+\eps-\la) } } {\Ga(n+1+\eps - ix) } \\
\times & \F{2}{1}{n+1+\eps+\la, n+\eps-\la}{n+1+\eps -
ix}{\frac{1}{1-e^{ -2i\phi}}}, \quad (n \in \Z).
\end{split}
\end{equation}
The parameters $\phi$ and $\eps$ satisfy the conditions $0<\phi<\pi$, $0\leq \eps <1$, and $\la$ satisfies one of the following conditions: $-\hf<\la<-\eps$,  $-\hf<\la<\eps-1$, or $\la=-\hf+i\rho$, $\rho \in \R$. In the last case $u_n$ is symmetric in $\rho$ and $-\rho$, so without loss of generality we assume $\rho \geq 0$. For $0<\phi\leq \frac{1}{6}\pi$ and $\frac{5}{6}\pi \leq \phi<\pi$ we use the unique analytic continuation of the $_2F_1$-function to $\C\setminus[1,\infty)$. Note that the Meixner-Pollaczek function is well defined for all $n \in \Z$, since $\Ga(c)^{-1} {}_2F_1(a,b;c;z)$ is analytic in $a$, $b$ and $c$.
 
The functions $u_n$ and $u_n^*$ satisfy the recurrence relation
\begin{equation} \label{rec MP funct}
2x \sin{\phi}\, u_n(x) = \al_n u_{n+1}(x)
 +  \be_n u_n(x) + \al_{n-1} u_{n-1}(x),
\end{equation}
where
\[
\al_n= \sqrt{(n+\eps+\lambda+1)(n+\eps- \lambda)}, \qquad
\be_n=2(n+\eps) \cos{\phi}.
\]
Let $L$ be the corresponding doubly infinite Jacobi operator acting on $\ell^2(\Z)$,
\[
L :e_n \mapsto \al_n e_{n+1} + \be_n e_n + \al_{n-1} e_{n-1}.
\]
$L$ is initially defined on finite linear combinations of the basis vectors of $\ell^2(\Z)$, and then $L$ is essentially self adjoint. The spectral measure of $L$ is described by
\[
\begin{split}
\inprod{u}{v} &= \frac{1}{2\pi} \int_\R \Big( \inprod{u}{U(x)}
\inprod{U(x)}{v} + \inprod{u}{U^*(x)} \inprod{U^*(x)}{v}\\ &- w_1(x)
\inprod{u}{U(x)} \inprod{U^*(x)}{v}- w_1^*(x)
\inprod{u}{U^*(x)} \inprod{U(x)}{v} \Big) w_0(x) dx,
\end{split}
\]
where
\begin{align*}
w_1(x) =& \frac{\Ga(\la+1+ix) \Ga(-\la+ix)} { \Ga(ix-\eps)
\Ga(1+\eps-ix)},\\ w_0(x) =& \frac{|\Ga(\la+1-ix) \Ga(-\la-ix)
\Ga(ix-\eps) \Ga(1+\eps-ix) |^2 }{ |\Ga(ix-\eps) \Ga(1+\eps-ix) |^2-
|\Ga(\la+1-ix) \Ga(-\la-ix)|^2 }
\\ & \times \frac{(2\sin \phi)^{-2\eps} e^{2x(\phi-\frac{\pi}{2})} }
{\Ga(-\eps-\la) \Ga(1+\la-\eps) \Ga(1+\la+\eps) \Ga(\eps-\la) }  \\
=&\, (2\sin \phi)^{-2\eps} e^{2x(\phi-\frac{\pi}{2})},\\
U(x) =& \sum_{n=-\infty}^\infty u_n(x;\la,\eps,\phi)e_n, \quad U^*(x) = \sum_{n=-\infty}^\infty u_n^*(x;\la,\eps,\phi)e_n.
\end{align*}
The spectral measure for $L$ can be obtained from \cite[\S4.4.11]{Koe}, using the connection formulas given there. The expression for $w_0(x)$ is found using Euler's reflection formula, elementary trigonometric identities and the conditions on $\la$.

Let $\mathcal H = \mathcal H(\la, \eps, \phi)$ be the Hilbert space consisting of functions 
\[
\mathbf{f}: \R \rightarrow \C^2,\quad
x \mapsto  \vect{f_1(x)}{f_2(x)},
\]
with inner product defined by
\[
\inprod{\mathbf{f}}{\mathbf{g}}=\frac{1}{2\pi} \int_\R \vect{g_1(x)}{g_2(x)}^*
\begin{pmatrix}
1 & -w_1(x) \\
-w_1^*(x)&1
\end{pmatrix}
\vect{f_1(x)}{f_2(x)} w_0(x) dx . 
\]
Observe that the square matrix inside the integral is positive definite and self adjoint. 
\begin{prop}\label{orthprop}
The functions
\[
\mathbf{u_n}=\left(\begin{matrix}u_n\\u_n^*\end{matrix}\right)
\]
form an orthonormal basis of the Hilbert space $\mathcal H$.
\end{prop}

The orthonormality follows by choosing $u$ and $v$ above
 as standard basis vectors $e_n$ and $e_m$. The completeness of the Meixner-Pollaczek functions follows from the uniqueness of the spectral measure. 
An alternative, group-theoretic approach, can be found in
\cite{VK}.  It is only worked out for the smaller range of parameters corresponding to $\mathrm{SU}(1,1)$ (rather than its Lie algebra or universal covering group). We will briefly describe how the  analytic arguments that lie behind the approach of \cite{VK} extend to the present situation. We only discuss the case $\lambda=-\frac 12+i\rho$, $\rho\in\mathbb R$, which is all that we 
need later.

Consider the space $L^2(\mathbb T)$ on the unit circle with respect to
normalized measure $|dz|/2\pi$. It has the orthonormal basis
$e_n(z)=z^n$, $n\in\mathbb Z$. 
If $\lambda=-\frac 12+i\rho$, $\rho\in\mathbb R$
and $\eps\in\mathbb R$, it is easily checked that
$$(Uf)(x)=\frac 1{\sqrt\pi}\,(x+i)^{\lambda-\eps}(x-i)^{\lambda+\eps}
f\left(\frac{x-i}{x+i}\right)$$
defines an isometry 
$U:\,L^2(\mathbb T)\rightarrow L^2(\mathbb R).$
Next we recall that the  Mellin transform, defined by
$$f(x)\mapsto\frac{1}{\sqrt{2\pi}} \int_0^\infty f(y) y^{ix-\frac 12}\,dy,$$
gives an  isometry 
$L^2(\mathbb R_{+})\rightarrow L^2(\mathbb R)$.
Thus, we may define a ``double'' Mellin transform as the isometry
$$V:L^2(\mathbb R)\rightarrow L^2(\mathbb R)\oplus L^2(\mathbb R)$$
given by
$$(Vf)(x)=\frac{1}{\sqrt{2\pi}}\left(\begin{matrix}
\int_0^\infty f(y)y^{ix-\frac 12}\,dy\\
\int_{0}^\infty f(-y)y^{ix-\frac 12}\,dy
\end{matrix}\right). $$

If we now let $T_t$, $t\in\mathbb R$,
 denote the translation operator $(T_tf)(x)=f(x+t)$, we
may compose the above isometries to obtain the orthonormal basis
\[
\mathbf{f_n}=\left(\begin{matrix}f_n^{+}\\f_n^{-}\end{matrix}\right)
=(V\circ T_t\circ U)(e_n)$$
 of  $L^2(\mathbb R)\oplus
L^2(\mathbb R)$. Explicitly, we have
$$f_n^{\pm}(x)=\frac{1}{\sqrt 2\,\pi}\int_0^\infty
(t+i\pm y)^{\lambda-n-\eps}
(t-i\pm y)^{\lambda+n+\eps}y^{ix-\frac 12}\,dy. 
\]
These  integrals may be expressed in terms of Gauss's hypergeometric
function; cf.~ \cite[\S 7.7.3]{VK}. Using Kummer's identities 
\cite[\S 2.9]{Erd}, one may then express $f_n^{\pm}$ in terms of
 $u_n(x+\rho)$ and $u_n^\ast(x+\rho)$, where $e^{2i\phi}=(t+i)/(t-i)$.
 Rewriting the identity $\langle
\mathbf{f_n},\mathbf{f_m}\rangle=\delta_{nm}$ in terms of Meixner-Pollaczek functions
and making a final change of variables $x\mapsto x-\rho$, one
recovers the orthonormal basis ${\mathbf u}_n$. 

The group-theoretic interpretation of this proof is the following. The
space $L^2(\mathbb T)$ is a natural representation space for the
principal unitary series ($e_n$ is proportional to the $e_n$ in
\eqref{princ} below). The operator $L$ gives the action of a 
hyperbolic Lie algebra element; cf.\ also \S4.2. It generates a
one-parameter subgroup of the universal covering group of $\mathrm
{SU}(1,1)$, which locally may be identified with the group of
linear fractional transformations of the circle that  have two common
fix-points. Thus, $L^2(\mathbb T)$  splits into two invariant
subspaces. The map $T_t\circ U$ corresponds to mapping the fix-points to
$\{0,\infty\}$, and the one-parameter subgroup to  dilations of
$\mathbb R$. Finally, the operator $V$ is the Fourier transform with
respect to these dilations.
In particular, the appearance of double eigenvalues in Proposition 
\ref{orthprop} has a natural geometric
explanation: it corresponds to the fact that a circle falls into
two pieces when removing two points.

\section{A bilinear summation formula} \label{sec2}
In this section we prove a bilinear summation formula for
Meixner-Pollaczeck polynomials. This summation is related to the
tensor product of a positive and a negative discrete series
representation of the Lie algebra $\su(1,1)$, which will be
explained in section \ref{sec:su11}. In the summation a certain
type of non-polynomial $_3F_2$-functions appear. These functions will be
investigated in section \ref{sec:Hahn}.

\begin{thm} \label{Thm:sum}
For $\rho^2 \in \mathrm{supp}\ d\mu(\cdot; k_2-k_1+\hf, k_1+k_2-\hf, k_1-k_2+p+\hf)$, $p \in \Z$, $x_1, x_2 \in \R$, and $k_1,k_2>0$, the Meixner-Pollaczek polynomials and the continuous dual Hahn polymials satisfy the following summation formula
\[
\begin{split}
\frac{1}{C}\sum_{n=0}^\infty & \frac{(-1)^n}{(2k_1)_{p+n} (2k_2)_n} s_n(\rho^2;k_2-k_1+\hf, k_1+k_2-\hf, k_1-k_2+p+\hf) p_{n+p}^{(k_1)}(x_1;\phi)
p_n^{(k_2)}(x_2;\phi) = \\
 & (-1)^p \frac{\Ga(i(x_1-x_2)+\hf+i\rho) \Ga(i(x_1-x_2)+\hf-i\rho) } {\Ga(k_2-k_1+i(x_1-x_2)+1) \Ga(k_1-k_2+i(x_1-x_2)+1+p) } \\
\times & \F{3}{2}{k_2-ix_2, k_2-k_1+\hf+i\rho, k_2-k_1+\hf-i\rho} {2k_2,
k_2-k_1+i(x_1-x_2)+1}{1} \\ \times & \F{2}{1}{p+k_1-k_2+\hf+i\rho,
p+k_1-k_2+\hf-i\rho}{ k_1-k_2+i(x_1-x_2)+p+1 } {
\frac{1}{1-e^{-2i\phi}}}\\
+ & \frac{ \Ga(i(x_2-x_1)+\hf+i\rho) \Ga(i(x_2-x_1)+\hf-i\rho) }
{\Ga(k_2-k_1+i(x_2-x_1)+1) \Ga(k_1-k_2+i(x_2-x_1)+1+p) } \\
\times & \F{3}{2}{k_2+ix_2, k_2-k_1+\hf+i\rho, k_2-k_1+\hf-i\rho}{2k_2,
k_2-k_1+i(x_2-x_1)+1}{1}\\
\times& \F{2}{1}{p+k_1-k_2+\hf+i\rho, p+k_1-k_2+\hf-i\rho}{ k_1-k_2+i(x_2-x_1)+p+1 }{
\frac{1}{1-e^{2i\phi}}},
\end{split}
\]
where
\[
C=\frac{e^{-2x_1(\phi-\frac{\pi}{2})}p!\, \Ga(2k_1)}{i^{p}\, (2\sin\phi)^{2k_1+p}\Ga(k_1+ix_1)
\Ga(k_1-ix_1)}.
\]
\end{thm}
Here we use the convention that if $p_n$ is a polynomial of degree $n$, then $p_{-n}=0$ for $n \in \N$. 
\begin{rem} \label{rem:sum}
(i) The $_2F_1$-series on the right hand side are Meixner-Pollaczek
functions, as defined by \eqref{def:MPfunctions}. The $_3F_2$-series on the right hand side can be considered as continuous Hahn functions, see section \ref{sec:Hahn}.

(ii) Theorem \ref{Thm:sum} is a non-polynomial extension of \cite[Thm.3.6]{GK}. Replacing $k_1+ix_1$, $k_2+ix_2$, $e^{2i\phi}$, $\rho^2$ in Theorem \ref{Thm:sum} by $x_1$, $x_2$, $c$, $y$ respectively, and restricting the new $x_1$ and $x_2$ to negative integers makes one term on the right hand side vanish, and we obtain \cite[Thm.3.6]{GK}.

(iii) There is an interesting limit case of the summation of Theorem \ref{Thm:sum}; for $x_1-x_2>0$
\[
\begin{split}
\sum_{n=0}^\infty &\frac{(-1)^n}{(2k_1)_{p+n} (2k_2)_n} s_n(\rho^2;k_2-k_1+\hf, k_1+k_2-\hf, k_1-k_2+p+\hf) L_{n+p}^{(2k_1-1)}(x_1) L_{n}^{(2k_2-1)}(x_2) = \\
&\Ga(2k_1) e^{x_2} \F{2}{1}{k_1+k_2-\hf+i\rho, k_1+k_2-\hf-i\rho}{2k_2}{\frac{x_2}{x_2-x_1}} \\
\times& (x_1-x_2)^{k_1+k_2-\hf-i\rho}U(p+k_1-k_2+\hf+i\rho;1+2i\rho;x_1-x_2), \end{split}
\]
where $L_n^{(\al)}$ is a Laguerre polynomial as defined in \cite{KS}, and $U(a;b;z)$ denotes the second solution of the confluent hypergeometric differential equation in the notation of Slater \cite{Sl}:
\[
U(a;b;z) = \frac{ \Ga(1-b)}{ \Ga(1+a-b) } \F{1}{1}{a}{b}{z} + \frac{\Ga(b-1)}{\Ga(a)} z^{1-b} \F{1}{1}{1+a-b}{2-b}{z}.
\]
This formula is obtained from Theorem \ref{Thm:sum} as follows. We replace $x_i$ by $-x_i/2\phi$, $i=1,2$, and transform the $_2F_1$-series on the right hand side by \cite[(2.3.12)]{AAR}. Then we let $\phi \rightarrow 0$. Here Stirling's formula is used, Euler's transformation is used for the $_2F_1$-series which are obtained from the $_3F_2$-series, and Kummer's transformation \cite[(1.4.1)]{Sl} is used for the $_1F_1$-series which are obtained from the $_2F_1$-series. This limit case can also be obtained by Lie algebraic methods, see \cite[Thm.3.10]{WG}.

(iv) Note that from  
\begin{equation} \label{p -p}
\begin{split}
s_n&(\rho^2;k_1-k_2+\hf, k_1+k_2-\hf, k_2-k_1-p+\hf) = \\
&(-1)^p |(k_1-k_2+\hf+i\rho)_p|^2 s_{n-p} (\rho^2;k_2-k_1+\hf, k_1+k_2-\hf, k_1-k_2+p+\hf),
\end{split}
\end{equation}
see \cite[(3.13)]{GK}, it follows that the sum on the left hand side of Theorem \ref{Thm:sum} is invariant under $k_1 \leftrightarrow k_2$, $x_1 \leftrightarrow x_2$, $p \leftrightarrow -p$.

(v) It is interesting to compare Theorem 3.1 with the results of \cite{IS}, where summation formulas with a similar, but simpler, structure are obtained for various orthogonal polynomials. The method used in \cite{IS} is completely different from the method we use here.
\end{rem}

\begin{proof}[Proof of Theorem \ref{Thm:sum}]
We start with the sum on the left hand side, with orthonormal polynomials:
\[
S=\sum_{n=0}^\infty S_n(\rho^2;k_2-k_1+\hf, k_1+k_2-\hf, k_1-k_2+p+\hf) P_{n+p}^{(k_1)}(x_1;\phi)
P_n^{(k_2)}(x_2;\phi).
\]
First we show that this sum converges absolutely. Writing out the summand $R_n$ explicitly gives
\[
R_n = K\,e^{-2in \phi} \frac{ \Ga( 2k_2+n) \Ga( 2k_1+p+n)} { \Ga(n+1) \Ga(n+p+1) }\ {}_3F_2 \ {}_2F_1 \ {}_2F_1
\]
where $K$ is a constant independent of $n$. To find the asymptotic behaviour for $n \rightarrow \infty$ of the $\Ga$-functions, we use the asymptotic formula for the ratio of two $\Ga$-functions \cite[\S4.5]{Olv}
\begin{equation}\label{assymp Ga1}
\frac{ \Ga (a+z)}{ \Ga(b+z) } = z^{a-b}\left( 1+
\mathcal O (z^{-1}) \right), \quad |z|\rightarrow \infty, \quad |\arg(z)|<\pi.
\end{equation}
The asymptotics for the $_3F_2$-function follows from transforming the function by \cite[p.15(2)]{Bai} and using \eqref{assymp Ga1}. This gives, for $n \rightarrow \infty$,
\[
\F{3}{2}{-n, k_1+k_2-\hf+i\rho, k_1+k_2-\hf-i\rho}{2k_2,2k_1+p}{1} = C_1\, n^{\hf-k_1-k_2-i\rho} + C_2\, n^{\hf-k_1-k_2+i\rho},
\] 
where $C_1$ and $C_2$ are independent of $n$. If $\rho^2$ is in the discrete part of $\mathrm{supp}\ d\mu(\cdot; k_2-k_1+\hf, k_1+k_2-\hf, k_1-k_2+p+\hf)$, we assume without loss of generality that $\Im(\rho)>0$. In this case the second term in the transformation \cite[p.15(2)]{Bai} vanishes, so $C_2=0$.  

The asymptotic behaviour of the $_2F_1$-functions follows from \cite[2.3.2(14)]{Erd}. This gives, for $0<\phi<\pi$ and $n\rightarrow \infty$,
\[
\begin{split}
\F{2}{1}{-(n+p), k_1 +ix_1}{2k_1}{1-e^{-2i\phi}} &= \left(C_3\, n^{-k_1-ix_1}+C_4\, n^{-k_1+ix_1} e^{-(n+p)(1-e^{-2i\phi})}\right) \left(1+\mathcal O(n^{-1}) \right),\\
\F{2}{1}{-n, k_2 -ix_2}{2k_2}{1-e^{2i\phi}} &= \left(C_5\, n^{-k_2-ix_2}+C_6\, n^{-k_2+ix_2} e^{-n(1-e^{-2i\phi})} \right)\left(1+\mathcal O(n^{-1}) \right),
\end{split}
\]
where $C_i$, $i=3,\ldots,6$, is independent of $n$. Since $\Re(1-e^{-2i\phi})=2\sin^2\phi>0$ and $x_1,x_2 \in \R$, we find
\[
R_n = (K_1\, n^{-\frac{3}{2}+i\rho} +K_2\, n^{-\frac{3}{2}-i\rho})\left(1+\mathcal O(n^{-1}) \right),\qquad n \rightarrow \infty,
\] 
and then for $\rho \in \R$ the sum $S$ converges absolutely. In case $\rho \in i\R$, we have $K_2=0$ and $\Im(\rho)>0$, and then $S$ still converges absolutely.\\

Next we write out the polynomials in $S$ as hypergeometric series, using \eqref{def:MP pol} and
\eqref{def:cont dHahn}, and then we transform the $_3F_2$-series, using the first formula on page 142 of \cite{AAR};
\begin{multline*}
\F{3}{2}{-n, k_2-k_1+\hf+i\rho, k_2-k_1+\hf-i\rho}{ 2k_2, p+1}{1} = \\
\frac{ (p+k_1-k_2+\hf-i\rho)_n } {(p+1)_n} \F{3}{2}{-n,
k_2-k_1+\hf+i\rho, k_1+k_2-\hf+i\rho} {2k_2, k_2-k_1+\hf+i\rho-n-p}{1}.
\end{multline*}
By \cite[\S2.9(27)]{Erd} with $(a,b,c,z)= (-n-p, k_1+ix_1, 2k_1, 1-e^{-2i\phi})$ the $_2F_1$-series for $P_{n+p}^{(k_1)}(x_1;\phi)$ is written as a sum of two $_2F_1$-series
\begin{equation} \label{eq:MPpol}
\begin{split}
&\F{2}{1}{ -n-p, k_1+ix_1}{ 2k_1} {1-e^{-2i\phi}}= \\&\quad
\frac{(1-e^{-2i\phi})^{-ix_1-k_1} \Ga(2k_1) \Ga(n+p+1) }{ \Ga(k_1-ix_1)
\Ga( n+p+1+k_1+ix_1) } \F{2}{1}{k_1+ix_1, 1-k_1+ix_1}{n+p+1+k_1+ix_1} {
\frac{1}{ 1-e^{-2i\phi}}}\\ &\quad+(e^{-2i\phi})^{n+p}
\frac{(1-e^{2i\phi})^{ix_1-k_1} \Ga(2k_1) \Ga(n+p+1) }{ \Ga(k_1+ix_1)
\Ga(n+p+1+k_1-ix_1) } \F{2}{1}{k_1-ix_1,
1-k_1-ix_1}{n+p+1+k_1-ix_1}{\frac{ 1}{1-e^{2i\phi}} }  .
\end{split}
\end{equation}
Now the sum $S$ splits according to this: $S=S_1+S_2$.

First we focus on $S_1$. Reversing the order of summation in the $_2F_1$-series of the second Meixner-Pollaczek polynomial and using Euler's transformation \cite[(2.2.7)]{AAR}, gives
\begin{multline} \label{step1}
\F{2}{1}{ -n, k_2+ix_2}{ 2k_2} {1-e^{-2i\phi}} = \\
 e^{-2in\phi}
(1-e^{2i\phi})^{ix_2-k_2} \frac{ (k_2+ix_2)_n } { (2k_2)_n }
\F{2}{1}{k_2-ix_2, 1-2k_2-ix_2}{1-k_2-n-ix_2}{\frac{1}{1-e^{-2i\phi}}}.
\end{multline}
Writing out the hypergeometric series as a sum, we get
\[
\begin{split}
S_1 =& e^{ip\phi} \frac{  (1-e^{-2i\phi})^{-k_1-ix_1}
(1-e^{2i\phi})^{ix_2-k_2} }{ \Ga(k_1-ix_1) \Ga(p+1+k_1+ix_1) }  \
\Ga(2k_1)\sqrt{ p! \, (2k_1)_p } \\ \times&
\sum_{n=0}^\infty \sum_{j=0}^n \sum_{l,m=0}^\infty  \frac{(k_2+ix_2)_n
(p+k_1-k_2+\hf-i\rho)_n (-n)_j (k_2-k_1+\hf+i\rho)_j
(k_1+k_2-\hf+i\rho)_j } { n!\,  (p+1+k_1+ix_1)_n \,j!\, (2k_2)_j
(k_2-k_1+\hf+i\rho-n-p)_j } \\ \times& \frac{ (k_1+ix_1)_l
(1-k_1+ix_1)_l (k_2-ix_2)_m (1-k_2-ix_2)_m } { l!\,
(n+p+1+k_1+ix_1)_l \, m! \, (1-k_2-n-ix_2)_m }
(1-e^{-2i\phi})^{-l-m}.
\end{split}
\]
Next we interchange summations
\[
\sum_{n=0}^\infty \sum_{j=0}^n = \sum_{j=0}^\infty \sum_{n=j}^\infty,
\]
then the sum over $n$ becomes
\[
\Sigma_1=\sum_{n=j}^\infty \frac{ (k_1-k_2+\hf-i\rho+p-j)_n
(k_2+ix_2-m)_n } { (n-j)! \, (p+1+l+k_1+ix_1)_n }.
\]
We substitute $k \mapsto n-j$, then we find by Stirling's formula for the summand $R$ of $S_1$, for large $j$ and $k$,  
\[
R \sim C \, j^{k_1-k_2+p} k^{-\hf} (j+k)^{-1-l-m},
\]
where $C$ is independent of $k$ and $j$. So we see that the sum $S_1$ converges absolutely for $k_1-k_2+p<0$, and $\frac{1}{6}\pi < \phi < \frac{5}{6}\pi$.
Now the sum $\Sigma_1$ is a multiple of a $_2F_1$-series with unit
argument, which is summable by Gauss' summation formula \cite[(46),
p.104]{Erd};
\[
\Sigma_1 = \frac{ (k_2-k_1+\hf-p+i\rho)_j (k_2+ix_2-m)_j
\Ga(p+l+1+k_1+ix_1) \Ga(l+m+\hf +ix_1-ix_2+i\rho) } {
 \Ga(l+j+\hf+k_2+ix_2+i\rho)
\Ga(p+l+m+k_1-k_2+ix_1-ix_2) }.
\]
This gives
\[
\begin{split}
S_1 = &e^{ip\phi} \frac{  (1-e^{-2i\phi})^{-k_1-ix_1}
(1-e^{2i\phi})^{ix_2-k_2} }{ \Ga(k_1-ix_1)  }  \ \Ga(2k_1)\sqrt{ p! \, (2k_1)_p } \\ \times & \sum_{l,m=0}^\infty \frac{
\Ga(l+m+i(x_1-x_2)+\hf+i\rho) (k_2-ix_2)_m (k_1+ix_1)_l (1-k_1+ix_1)_l
} { \Ga(p+l+m+1+k_1-k_2+i(x_1-x_2)) \Ga(l+k_2+\hf+ix_1+i\rho) \, m! \,
l! } (1-e^{-2i\phi})^{-l-m} \\ \times & \sum_{j=0}^\infty \frac{
(k_2-k_1+\hf+i\rho)_j (k_1+k_2-\hf+i\rho)_j (k_2+ix_2-m)_j } {j!\,
(2k_2)_j (l+\hf+k_2+ix_1+i\rho)_j }
\end{split}
\]
The sum over $j$ is a $_3F_2$-series, which by Kummer's transformation \cite[Cor.3.3.5]{AAR} becomes
\[
\begin{split}
&\frac{ \Ga(l+\hf+k_2+ix_1+i\rho)
\Ga(l+m+\hf+i(x_1-x_2)-i\rho) }{ \Ga(k_2-k_1+1+i(x_1-x_2)+l+m)
\Ga(l+k_1+ix_1) } \\& \qquad \times \F{3}{2}{k_2-k_1+\hf+i\rho,
k_2-k_1+\hf-i\rho, k_2-ix_2+m} {2k_2, k_2-k_1+1+i(x_1-x_2)+l+m}{1}.
\end{split}
\]
Here we need the condition $k_1>0$ for absolute convergence.
We write out the $_3F_2$-series explicitly as a sum over $j$ and
interchange summations
\[
S_1 = \sum_{l,m,j=0}^\infty = \sum_{n=0}^\infty \sum_{j=0}^\infty
\sum_{l+m=n}.
\]
For the sum over $l+m=n$ we find
\[
\sum_{m+l=n} \frac{ (1-k_1+ix_1)_l (k_2-ix_2+j)_m } {l!\ m!} = \frac{ (
1+k_2-k_1+i(x_1-x_2)+j)_n }{n!}.
\]
Now $S_1$ reduces to a double sum, which splits as a product of two sums, and we obtain
\[
\begin{split}
S_1 = &
e^{ip\phi} \frac{ (1-e^{-2i\phi})^{-k_1-ix_1} (1-e^{2i\phi})^{ix_2-k_2}
\Ga(\hf+i(x_1-x_2)+i\rho) \Ga(\hf+i(x_1-x_2)-i\rho)}{ \Ga(k_1-ix_1)
\Ga(k_1+ix_1) \Ga(p+k_1-k_2+i(x_1-x_2)+1) \Ga(k_2-k_1+i(x_1-x_2)+1)} \\
& \times \Ga(2k_1)\sqrt{ p! \, (2k_1)_p }
\F{2}{1}{\hf+i(x_1-x_2)+i\rho, \hf+i(x_1-x_2)-i\rho}
{p+k_1-k_2+i(x_1-x_2)+1} {\frac{1}{1-e^{-2i\phi}}} \\ & \times
\F{3}{2}{k_2-k_1+\hf+i\rho, k_2-k_1+\hf-i\rho, k_2-ix_2} {2k_2,
k_2-k_1+i(x_1-x_2)+1}{1}.
\end{split}
\]

$S_2$ is calculated in the same way as $S_1$, only for the first step  \cite[\S2.10(4)]{Erd} is used instead of \eqref{step1}. Then we obtain
\[
\begin{split}
S_2 &=  e^{-ip\phi}\frac{(1-e^{-2i\phi})^{-k_2-ix_2} (1-e^{2i\phi})^{ix_1-k_1} } {\Ga(k_1+ix_1) \Ga(k_1-ix_1)} \Ga(2k_1)\sqrt{ p! \, (2k_1)_p }\\
&\times \frac{ \Ga(i(x_2-x_1)+\hf+i\rho) \Ga(i(x_2-x_1)+\hf-i\rho) }
{\Ga(k_2-k_1+i(x_2-x_1)+1) \Ga(k_1-k_2+i(x_2-x_1)+1+p) } \\  &\times
\F{3}{2}{k_2+ix_2, k_2-k_1+\hf+i\rho, k_2-k_1+\hf-i\rho}{2k_2,
k_2-k_1+i(x_2-x_1)+1}{1}\\
 & \times\F{2}{1}{\hf+i(x_2-x_1)+i\rho,
\hf+i(x_2-x_1)-i\rho}{ k_1-k_2+i(x_2-x_1)+p+1 }{ \frac{1}{1-e^{2i\phi}}}.
\end{split}
\]
Finally using Euler's transformation for the $_2F_1$-series, using
\begin{align*}
e^{ip\phi} (1-e^{2i\phi})^{-p} &= (-1)^p(2i \sin\phi)^{-p}, \\
(1-e^{-2i\phi})^{-ix_1}(1-e^{2i\phi})^{ix_1} &=  e^{-2x_1(\phi-\frac{\pi}{2})},\\
(1-e^{-2i\phi})^{-k_1}(1-e^{2i\phi})^{-k_1} &=  (2\sin\phi)^{-2k_1},
\end{align*}
and writing the polynomials in the normalization given by \eqref{def:MP pol} and \eqref{def:cont dHahn}, the theorem is proved.

Let us remark that all the series used in the proof are absolutely convergent under the conditions $\rho^2 \in \mathrm{supp}\ d\mu(\cdot; k_2-k_1+\hf, k_1+k_2-\hf, k_1-k_2+p+\hf)$, $x_1, x_2 \in \R$, $k_1>0$, $\frac{1}{6}\pi<\phi<\frac{5}{6}\pi$ and $k_1-k_2+p<0$. The last condition can be removed using the symmetry $(k_1,k_2,p) \leftrightarrow (k_2, k_1, -p)$ and continuity in $k_1$ and $k_2$. Using the analytic continuation of the hypergeometric function, we see that the result remains valid for $0 < \phi < \pi$.
\end{proof}

\section{Clebsch-Gordan coefficients for hyperbolic basis vectors of $\su(1,1)$} \label{sec:su11}
In this  section we consider the tensor product of a positive and a negative discrete series representation. We diagonalize a certain self-adjoint element of $\su(1,1)$ using (doubly infinite) Jacobi operators. We also give generalized eigenvectors, which can be considered as hyperbolic basis vectors. Using the summation formula from the previous section, we show that the Clebsch-Gordan coefficients for the eigenvectors are continuous Hahn functions. We find the corresponding integral transform pair by formal computations. In order to give a rigorous proof for the continuous Hahn integral transform, we realize the generators of $\su(1,1)$ in the discrete series as difference operators acting on polynomials. Using these realizations, the Casimir element in the tensor product is realized as a difference operator. Spectral analysis of this difference operator is carried out in section \ref{sec:Hahn}.

\subsection{The Lie algebra $\boldsymbol{\su(1,1)}$} \label{ssec:su11}

The Lie algebra $\su(1,1)$ is generated by the elements $H$, $B$ and $C$, satisfying the commutation relations
\begin{equation} \label{commutation rel}
[H,B]=2B, \quad [H,C]=-2C, \quad [B,C]=H.
\end{equation}
There is a $*$-structure defined by $H^*=H$ and $B^*=-C$. 
The center of $\mathcal U\big(\su(1,1)\big)$ is generated by the Casimir element $\Omega$, which is given by 
\begin{equation} \label{Casimir}
 \Omega = -\frac{1}{4}(H^2+2H+4CB)=-\frac{1}{4}(H^2-2H+4BC).
\end{equation}

There are four classes of irreducible unitary representations of
$\su(1,1)$, see \cite[\S6.4]{VK}:

The positive discrete series
representations $\pi_k^+$ are representations labelled by $k>0$. The
representation space is $\lt(\Z_{\geq 0})$ with orthonormal basis
$\{e_n\}_{n\in \Z_{\geq 0}}$. The action is given by
\begin{equation} \label{pos}
\begin{split}
\pi_k^+(H)\, e_n =&\ 2(k+n)\, e_n,  \\ \pi_k^+(B)\, e_n =&\
\sqrt{(n+1)(2k+n)}\, e_{n+1},  \\ \pi_k^+(C)\, e_n =&\ -\sqrt{n(2k+n-1)}\,
e_{n-1},  \\ \pi_k^+(\Omega)\,e_n =&\ k(1-k)\, e_n.
\end{split}
\end{equation}

The negative discrete series representations $\pi_{k}^-$ are labelled by
$k>0$. The representation space is $\lt(\Z_{\geq 0})$ with orthonormal
basis $\{e_n\}_{n \in \Z_{\geq 0}}$. The action is given by
\begin{equation} \label{neg}
\begin{split}
\pi_{k}^-(H)\, e_n =&\ -2(k+n)\,e_n, \\ \pi_{k}^-(B)\, e_n =&\
-\sqrt{n(2k+n-1)}\, e_{n-1},  \\ \pi_{k}^-(C)\, e_n =&\ \sqrt{(n+1)(2k+n)}\,
e_{n+1},  \\ \pi_{k}^-(\Omega)\, e_n =&\ k(1-k)\, e_n.
\end{split}
\end{equation}

The principal series representations $\pi^{\rho,\eps}$ are labelled by
$\eps \in [0,1)$ and $\rho \geq 0$, where $(\rho,\eps) \neq
(0,\frac{1}{2})$. The representation space is $\lt(\Z)$ with
orthonormal basis $\{e_n\}_{n \in \Z}$. The action is given by
\begin{equation} \label{princ}
\begin{split}
\pi^{\rho, \eps}(H)\, e_n =&\ 2(\eps+n)\, e_n, \\ \pi^{\rho, \eps}(B)\, e_n
=&\ \sqrt{(n+\eps+\frac{1}{2}-i\rho) (n+\eps+\frac{1}{2}+i\rho)}\,
e_{n+1},
\\ \pi^{\rho, \eps}(C)\, e_n =&\
-\sqrt{(n+\eps-\frac{1}{2}-i\rho) (n+\eps-\frac{1}{2}+i\rho)}\, e_{n-1},
\\ \pi^{\rho, \eps}(\Omega)\, e_n =&\ (\rho^2+\frac{1}{4})\, e_n.
\end{split}
\end{equation}
For $(\rho,\eps)= (0,\hf)$ the representation $\pi^{0,\hf}$ splits into a direct sum of a positive and a negative discrete series representation: $\pi^{0,\hf}= \pi^+_{1/2} \oplus \pi^-_{1/2}$. The representation space splits into two invariant subspaces: $\{e_n \, |\, n<0 \} \oplus \{ e_n\, | \, n \geq 0 \}$.

The complementary series representations $\pi^{\lambda,\eps}$ are
labelled by $\eps$ and $\lambda$, where $\eps \in [0,\frac{1}{2})$ and
$\lambda \in (-\frac{1}{2},-\eps)$  or $\eps \in (\frac{1}{2},1)$ and
$\lambda \in (-\frac{1}{2},\eps-1)$. The representation space is
$\lt(\Z)$ with orthonormal basis $\{e_n\}_{n \in \Z}$. The action is
given by
\begin{equation} \label{comp}
\begin{split}
\pi^{\lambda,\eps}(H)\, e_n =&\ 2(\eps+n)\, e_n, \\ \pi^{\lambda,\eps}(B)\,
e_n =&\ \sqrt{(n+\eps+1+\lambda) (n+\eps-\lambda)}\, e_{n+1}, \\
\pi^{\lambda,\eps}(C)\, e_n =&\ -\sqrt{(n+\eps+\lambda)
(n+\eps-\lambda-1)}\, e_{n-1},  \\ \pi^{\lambda,\eps}(\Omega)\, e_n =&\
-\lambda(1+\lambda)\, e_n.
\end{split}
\end{equation}
Note that if we formally write $\la = -\frac{1}{2}+i\rho$ the actions in the principal series and in the complementary series are the same.

We remark that the operators \eqref{pos}-\eqref{comp} are unbounded, with domain the set of finite linear combinations of the basis vectors. The representations are $*$-representations in the sense of Schm\"udgen \cite[Ch.8]{Sch}.\\

The decomposition of the tensor product of a positive and a
negative discrete series representation of $\su(1,1)$ is determined in full generality in \cite[Thm.2.2]{GK}.
\begin{thm} \label{thm:decomp}
For $k_1 \leq k_2$ the decomposition of the tensor product of positive
and negative discrete series representations of $\su(1,1)$ is
\begin{align*}
\pi_{k_1}^+ \tensor\pi_{k_2}^- &\cong \dirint \pi^{\rho,\eps} d \rho, &
k_1-k_2 \geq -\frac{1}{2}, k_1+k_2 \geq \frac{1}{2}, \\ \pi_{k_1}^+
\tensor\pi_{k_2}^- &\cong \dirint \pi^{\rho,\eps} d \rho \oplus
\pi^{\la, \eps}, & k_1+k_2<\frac{1}{2}, \\ \pi_{k_1}^+
\tensor\pi_{k_2}^- &\cong \dirint \pi^{\rho,\eps} d \rho \oplus
\bigoplus_{\substack{j \in \Z_{\geq 0}\\ k_2-k_1-\frac{1}{2}-j>0}}
\pi_{k_2-k_1-j}^-, & k_1-k_2
< -\frac{1}{2},
\end{align*}
where $\eps = k_1-k_2+L$, $L$ is the unique integer such that $\eps \in
[0,1)$, and $\la = -k_1-k_2$. 
The intertwiner $J$ is given by
\begin{equation} \label{decomp vec}
J(e_{n_1}\tensor e_{n_2}) =(-1)^{n_2}\int_\R S_n(y;n_1-n_2)\, e_{n_1-n_2-L}\, d\mu^\hf(y;n_1-n_2),
\end{equation}
where $n= \min\{n_1,n_2\}$, $S_n(y;p)$ is an orthonormal continuous dual Hahn polynomial,
\[
S_n(y;p) =
\begin{cases}
S_n(y;k_1-k_2+\frac{1}{2},k_1+k_2-\frac{1}{2},k_2-k_1-p+\frac{1}{2}),&
p \leq 0,\\
S_n(y;k_2-k_1+\frac{1}{2},k_1+k_2-\frac{1}{2},k_1-k_2+p+\frac{1}{2}),&
p\geq0,
\end{cases}
\]
and $d\mu(y;p)$ is the corresponding orthogonality measure
\[
d \mu(y;p) =
\begin{cases}
d\mu(y;k_1-k_2+\frac{1}{2},k_1+k_2-\frac{1}{2},k_2-k_1-p+\frac{1}{2}),&
p \leq 0,\\ d \mu
(y;k_2-k_1+\frac{1}{2},k_1+k_2-\frac{1}{2},k_1-k_2+p+\frac{1}{2}),&
p\geq0.
\end{cases}
\]
\end{thm}
The inversion of \eqref{decomp vec} can be given explicitly, e.g.~ for an element
\[
f \tensor e_{r-L} = \int_0^\infty f(x) e_{r-L} dx \in L^2(0,\infty)
\tensor \lt(\Z) \cong \dirint \lt(\Z) dx
\]
in the representation space of the direct integral representation, we have
\begin{equation} \label{inv decomp}
J^*(f \tensor e_{r-L}) =
\begin{cases}
\displaystyle \sum_{p=0}^\infty (-1)^{p-r} \left[ \int_{\R} S_p(y;r)
f(y) d\mu^\hf(y;r) \right] e_p \tensor e_{p-r}, & r \leq 0,\\ \displaystyle
\sum_{p=0}^\infty (-1)^{p} \left[\int_{\R} S_p(y;r) f(y) d \mu^\hf(y;r)
\right] e_{p+r} \tensor e_p, & r \geq 0.
\end{cases}
\end{equation}
For the discrete components in Theorem \ref{thm:decomp} we can replace $f$ by a Dirac delta function at the appropriate points of the discrete mass of $d \mu(\cdot;r)$. We remark that for $k_1=k_2<1/4$, the occurrence of a complementary series representation in the tensor product was discovered
by Neretin \cite{Ne}. This phenomenon was investigated from the viewpoint of operator theory in \cite{EHPRZ}.

In the following subsection we assume that discrete terms do not occur in the tensor product decomposition. From the calculations it is clear how to extend the results to the general case. At the end of section \ref{sec:Hahn} we briefly discuss the results for the discrete terms in the decomposition. \\

In the Lie algebra $\su(1,1)$ three types of elements can be distinguished: the elliptic, the parabolic and the hyperbolic elements. These are related to the three conjugacy classes of the group $SU(1,1)$. A basis on which an elliptic element acts diagonally is called an elliptic basis, and similarly for the parabolic and hyperbolic elements. The basisvectors $e_n$ in \eqref{pos}-\eqref{comp} are elliptic basisvectors.

We consider self-adjoint elements of the form 
\[
-aH + B - C \in \su(1,1), \qquad (a \in \R),
\]
in the tensor product of a positive and a negative discrete series representation. For $|a|=1$ this is a parabolic element, for $|a|<1$ it is hyperbolic and for $|a|>1$ it is elliptic. For the elliptic and the parabolic case we refer to \cite{GK}, respectively \cite{WG}. We consider the case $|a|<1$.

\subsection{Hyperbolic basisvectors} \label{ssec:H}

We consider a self-adjoint element ${X_\phi}$ in
$\mathfrak{su}(1,1)$, given by
\[
X_\phi = -\cos{\phi}\ H + B -C, \qquad 0<\phi<\pi.
\]
The action of $X_\phi$ in the discrete series can be identified with the three-term recurrence relation for the Meixner-Pollaczek polynomials \eqref{MP rec}, cf. \cite[Prop.~3.1]{KJ1}.
\begin{prop}\label{prop:+-H}
The operators $\Te^\pm$, defined by
\[
\begin{split}
\Te^\pm : \ell^2(\Z_{\geq 0}) &\rightarrow L^2\big(\R, w^{(k)}(x;\phi)dx\big)\\
e_n &\mapsto P_n^{(k)}(\cdot;\phi),
\end{split}
\]
are unitary and intertwine $\pi^\pm_k(X_\phi)$ with $M(\pm 2x \sin \phi)$.
\end{prop}
Here $M$ denotes the multiplication operator, i.e. $M(f)g(x)=f(x)g(x)$.

Proposition \ref{prop:+-H} states that
\[
v^\pm(x) = \sum_{n=0}^\infty P_n^{(k)}(x;\phi)\, e_n,
\]
are generalized eigenvectors of $\pi^\pm_k(X_\phi)$ for eigenvalue $\pm 2x \sin{\phi}$. These eigenvectors can be considered as hyperbolic basis vectors.

The action of $X_\phi$ in the principal unitary series can be identified with the recurrence relation for the Meixner-Pollaczek functions \eqref{rec MP funct}. Then the spectral decomposition of the corresponding doubly infinite Jacobi operator gives the following.
\begin{prop}\label{prop:princ H}
The operator $\Te^{\rho,\eps}$ defined by
\[
\begin{split}
\Te^{\rho,\eps} : \ell^2(\Z) &\rightarrow \mathcal H(-\hf+i\rho, \eps, \pi-\phi), \\
e_n &\mapsto \vect{ u_n(\cdot;-\hf+i\rho, \eps, \pi-\phi) }{  u^*_n(\cdot;-\hf+i\rho, \eps,\pi- \phi)},
\end{split}
\]
is unitary and intertwines $\pi^{\rho,\eps}(X_\phi)$ with $M(2x\sin\phi)$, and extends to a unitary equivalence.
\end{prop}
From Propostion \ref{prop:princ H} we obtain that
\[
\vect{v_{\rho,\eps}(x)}{v_{\rho,\eps}^*(x)} =  \sum_{n=-\infty}^\infty \vect{ u_n(x;-\hf+i\rho, \eps, \pi-\phi) }{  u^*_n(x;-\hf+i\rho, \eps,\pi- \phi)} \, e_n
\]
is a generalized eigenvector of $\pi^{\rho,\eps}(X_\phi)$ for eigenvalue $2x\sin\phi$. 

Next we consider the action of $X_\phi$ in the tensor product.
Recall that in the tensor product we need the coproduct $\De$, defined by $\De(Y) = 1 \tensor Y + Y \tensor 1$ for $Y \in \su(1,1)$. Then from Proposition \ref{prop:+-H} we find the following.
\begin{prop} \label{prop:+-1}
The operator $\Ups$ defined by
\[
\begin{split}
\Ups : \ell^2(\Z_{\geq 0}) \tensor \ell^2(\Z_{\geq 0}) &\rightarrow L^2\big(\R^2, w^{(k_1)}(x_1;\phi)w^{(k_2)}(x_2;\phi)dx_1 dx_2\big),\\
e_{n_1} \tensor e_{n_2} &\mapsto P_{n_1}^{(k_1)}(x_1;\phi) P_{n_2}^{(k_2)}(x_2;\phi),
\end{split}
\]
is unitary and intertwines $\pi^+_{k_1} \tensor \pi^-_{k_2} \big( \De(X_\phi) \big)$ with $M\big(2(x_1-x_2)\sin\phi\big)$.
\end{prop}
In terms of the generalized eigenvectors $v^+$ and $v^-$, we find from Proposition \ref{prop:+-1} that
\[
v^+(x_1) \tensor v^-(x_2)=\sum_{n_1,n_2=0}^\infty P_{n_1}^{(k_1)}(x_1;\phi)
P_{n_2}^{(k_2)}(x_2;\phi)\, e_{n_1} \tensor e_{n_2},
\]
is a generalized eigenvector of $\pi^+_{k_1} \tensor \pi^-_{k_2}\big(\De(X_\phi)\big)$ for eigenvalue $2(x_1-x_2)\sin\phi$.\\

To determine the action of $\Ups$ on the representation space of the direct integral representation $\int^\oplus \pi^{\rho,\eps} d\rho$, we need to find the operator $\tilde \Ups$, such that $\Ups = \tilde \Ups \circ J$. Here $J$ is the intertwiner defined in Theorem \ref{thm:decomp}. For appropriate functions $g_1$ and $g_2$ we define an operator $\tilde \Ups_g$ by
\[
\begin{split}
\tilde \Ups_{\mathbf{g}} : L^2(0,\infty)
\tensor \lt(\Z) \cong \dirint \lt(\Z) dx &\rightarrow \dirint \mathcal H(-\hf+i\rho,\eps, \pi-\phi)\, d\rho,\\
f \tensor e_n & \mapsto \int_0^\infty f(\rho) \vect{g_1(\rho)}{g_2(\rho)}^* \vect{ u_n(t;-\hf+i\rho, \eps, \pi-\phi) }{  u^*_n(t;-\hf+i\rho, \eps,\pi- \phi)} d\rho.
\end{split}
\]
From Proposition \ref{prop:princ H} we see that $\tilde \Ups_{\mathbf{g}}$ intertwines $\int^\oplus \pi^{\rho,\eps}(X_\phi) d\rho$ with $M\big(2t\sin\phi\big)$. The functions $g_1$ and $g_2$ for which $\Ups = \tilde \Ups_{\mathbf{g}} \circ J$ are the Clebsch-Gordan coefficients for the hyperbolic bases. To determine the Clebsch-Gordan coefficients we use the summation formula in Theorem \ref{Thm:sum}. Define the continuous Hahn function by
\begin{equation} \label{cont.H.f}
\begin{split}
\ph_\rho(x;t)= &\,\ph_\rho(x;t,k_1,k_2,\phi)\\
 =&\, \frac{ e^{-x(2\phi-\pi)}} { |\Ga(k_1+ix+it)|^2} \F{3}{2}{k_2-ix, k_2-k_1+\hf+i\rho,
k_2-k_1+\hf-i\rho} {2k_2, k_2-k_1+it+1} {1}.
\end{split}
\end{equation}
\begin{thm} \label{thm:CG-decomp H}
Let
\[
\vect{g_1(\rho)}{g_2(\rho)}^* = 
\vect{\ph_\rho(x_2;x_1-x_2)}{\ph_\rho^*(x_2;x_1-x_2)}^*
\begin{pmatrix}
m(\rho) & 0 \\
0 & \overline{m(\rho)}
\end{pmatrix},
\]
where $\ph_\rho(x_2;x_1-x_2)=\ph_\rho(x_2;x_1-x_2,k_1,k_2,\phi)$ and $m(\rho)$ is given by
\[
\begin{split}
m(\rho) =& (-i)^L\frac{(2\sin\phi)^{-2k_1-L}e^{(x_2-x_1)(2\phi-\pi)}}{\sqrt{2\pi}} \frac{ \Ga(\hf+i(x_2-x_1)+i\rho)
\Ga(\hf+i(x_2-x_1)-i\rho) } {\Ga(k_2-k_1+i(x_2-x_1)+1) }\\
&\times \sqrt{\frac{\Ga(2k_1)}{\Ga(2k_2)}}  \left| \frac{ \Ga(k_1+k_2-\hf+i\rho)
\Ga(k_2-k_1+\hf+i\rho) } {\Ga(2i\rho)} \right|,
\end{split}
\]
then we have $\Ups = \tilde \Ups_{\mathbf{g}} \circ J$.
\end{thm}
\begin{proof}
To show that $\Ups = \tilde \Ups_{\mathbf{g}} \circ J$,
we use the summation formula of Theorem \ref{Thm:sum} with the orthogonal polynomials written in orthonormal form. We multiply by a continuous dual Hahn polynomial of degree $n_2$ with the same parameters as in Theorem \ref{Thm:sum}, with $p=n_1-n_2 \geq 0$.  Then integrating against the corresponding orthogonality measure gives an equality with the following structure
\[
P_{n_1}^{(k_1)}(x_1;\phi) P_{n_2}^{(k_2)}(x_2;\phi) = \int_0^\infty S_{n_2}(\rho) ( {}_3F_2\ {}_2F_1 + {}_3F_2\ {}_2F_1) d\rho.
\]
The $_2F_1$-functions are the Meixner-Pollaczek functions as defined by \eqref{def:MPfunctions}.
From Proposition \ref{prop:+-H} we see that the left hand side is equal to $\Ups(e_{n_1} \tensor e_{n_2})(x_1,x_2)$. So from Theorem \ref{thm:decomp}, with $n_1 \geq n_2$, and Proposition \ref{prop:princ H} it follows that the right hand side must be equal to 
\[
 \int_0^\infty S_{n_2}(\rho) \mathbf{g}^*(\rho)    \Te^{\rho,\eps}(e_{n_1-n_2-L})(x_1-x_2) d\mu^\hf(\rho) ,
\]
where $\mathbf{g}$ is a vector containing the Clebsch-Gordan coefficients for the hyperbolic bases. This gives the desired result. For $n_1-n_2<0$ the theorem follows after using \eqref{p -p}.
\end{proof}

\begin{rem}
The explicit expressions of the Clebsch-Gordan coefficients as $_3F_2$-series can also be found in Mukunda and Radhakrishnan \cite{MuR}. The method used in \cite{MuR} is completely different from the method used here.
\end{rem}

In terms of the generalized eigenvectors, Theorem \ref{thm:CG-decomp H} states
\[
J\big(v^+(x_1) \tensor v^-(x_2) \big)= \int_0^\infty \vect{\ph_\rho(x_2;x_1-x_2)}{\ph_\rho^*(x_2;x_1-x_2)}^*
\begin{pmatrix}
m(\rho) & 0 \\
0 & \overline{m(\rho)}
\end{pmatrix}
\vect{v_{\rho,\eps}(x_1-x_2)}{v_{\rho, \eps}^*(x_1-x_2)} d\rho.
\]

\subsection{The continuous Hahn integral transform} \label{ssec:formal}
Since the continuous Hahn functions occur as Clebsch-Gordan coefficients for hyperbolic bases, they should satisfy (generalized) orthogonality relations. We find these relations by formal computations with the generalized eigenvectors. \\

For an element $f \in \ell^2(\Z_{\geq 0}) \tensor \ell^2(\Z_{\geq 0})$ we have the transform pair
\begin{equation} \label{transf1}
\begin{cases}
\displaystyle (\Ups f)(x_1,x_2) = \inprod{f}{v^+(x_1) \tensor v^-(x_2) },\\ \\
\displaystyle f = \iint_{\R^2}  (\Ups f)(x_1,x_2)\ v^+(x_1) \tensor v^-(x_2)\ w^{(k_1)}(x_1;\phi) w^{(k_2)}(x_2;\phi)dx_1 dx_2.
\end{cases}
\end{equation}
Similarly for $f  \in \ell^2(\Z)$ we have the transform pair
\begin{equation} \label{transf2}
\begin{cases}
\displaystyle (\Te^{\rho,\eps} f)(x)=\left\langle f , \vect{v_{\rho,\eps}(x)} {v_{\rho,\eps}^*(x)} \right\rangle \in \C^2,\\ \\
\displaystyle f=\frac{1}{2\pi} \int_\R (\Te^{\rho,\eps} f) (x)^*
\begin{pmatrix}
1 & -w_1(x) \\
-w_1^*(x)&1
\end{pmatrix}
\vect{v_{\rho,\eps}(x)} {v_{\rho,\eps}^*(x)} w_0(x) dx . 
\end{cases}
\end{equation}
Denoting the intertwiner $J$ in Theorem \ref{thm:decomp} by $J= \dirint J_\rho d\rho$, we find from Theorem \ref{thm:CG-decomp H}
\[
(\Ups f)(x_1,x_2) = \int_0^\infty \vect{\ph_\rho(x_2;x_1-x_2)}{\ph_\rho^*(x_2;x_1-x_2)}^*
\begin{pmatrix}
m(\rho) & 0 \\
0 & \overline{m(\rho)}
\end{pmatrix}
\big(\Te^{\rho,\eps}(J_\rho f)\big)(x_1-x_2) d\rho.
\]
We want to invert this formula. From \eqref{transf1} and Theorem \ref{thm:CG-decomp H} we find formally
\[
\begin{split}
J_\rho f = \iint_{\R^2} &(\Ups f)(x_1,x_2) w^{(k_1)}(x_1;\phi) w^{(k_2)}(x_2;\phi)\\ \times &\vect{\ph_\rho(x_2;x_1-x_2)}{\ph_\rho^*(x_2;x_1-x_2)}^*
\begin{pmatrix}
m(\rho) & 0 \\
0 & \overline{m(\rho)}
\end{pmatrix} 
\vect{v_{\rho,\eps}(x_1-x_2)} {v_{\rho,\eps}^*(x_1-x_2)}
dx_1 dx_2.
\end{split}
\]
We substitute $x_1 \mapsto x+t$, $x_2 \mapsto x$, then \eqref{transf2} with $f$ replaced by $J_\rho f$ gives
\[
\begin{split}
w_0(t) &\begin{pmatrix}
1 & -w_1^*(t) \\
-w_1(t)&1
\end{pmatrix}
\big(\Te^{\rho,\eps} (J_\rho f)\big)(t) = \\
&2\pi \int_\R (\Ups f)(x+t,x) 
\begin{pmatrix}
\overline{m(\rho)} & 0 \\
0 & m(\rho)
\end{pmatrix}
\vect{\ph_\rho(x;t)}{\ph_\rho^*(x;t)}
w^{(k_1)}(x+t;\phi) w^{(k_2)}(x;\phi)dx.
\end{split}
\]
We denote  
\[
(\Te^{\rho,\eps} \big(J_\rho f)\big)(t) = \frac{2\pi}{w_0(t)} 
\begin{pmatrix}
\overline{m(\rho)} & 0 \\
0 & m(\rho)
\end{pmatrix}
\begin{pmatrix}
1 & -w_1^*(t) \\
-w_1(t)&1
\end{pmatrix}^{-1}
\mathbf{g}(\rho)
\]
and $(\mathcal F \mathbf{g})(x) =(\Ups f)(x+t,x)$, then we have the following integral transform pair with the continuous Hahn functions as a kernel
\[
\begin{cases}
\displaystyle (\mathcal F \mathbf{g})(x)= \frac{2\pi}{w_0(t)} \int_0^\infty \vect{\ph_\rho(x;t)}{\ph_\rho^*(x;t)}^*
\begin{pmatrix}
1 &  w_1^*(t) \\
w_1(t)& 1
\end{pmatrix}
\mathbf{g}(\rho)\, \frac{|m(\rho)|^2}{1-|w_1(t)|^2}d\rho\\ \\
\displaystyle \mathbf{g}(\rho)=
 \int_\R (\mathcal F \mathbf{g})(x)
\vect{\ph_\rho^*(x;t)}{\ph_\rho(x;t)}
w^{(k_1)}(x+t;\phi) w^{(k_2)}(x;\phi)dx.
\end{cases}
\]

In section \ref{sec:Hahn} we give a rigorous proof for this integral transform pair using spectral analysis of a difference operator for which the continuous Hahn functions are eigenfunctions. In the next subsection we obtain the difference operator from the action of the Casimir element in the tensor product.

\subsection{A realization of the discrete series representations.} 
The following lemma is based on the fact that $\mathfrak{sl}(2,\C)$ is semi-simple, so $[\mathfrak{sl}(2,\C),\mathfrak{sl}(2,\C)] = \mathfrak{sl}(2,\C).$
\begin{lem} \label{lem:HXphi}
\[
B = \frac{1}{4}[H,X_\phi]+\hf X_\phi+\hf\cos\phi \, H, \quad
C=\frac{1}{4}[H,X_\phi]-\hf X_\phi-\hf\cos\phi \,H.
\]
\end{lem}
\begin{proof}
This follows from the definition of $X_\phi$ and the commutation relations \eqref{commutation rel}.
\end{proof}
This lemma shows that to find the action of the generators $H$, $B$ and $C$ on the Meixner-Pollaczek polynomials, it is enough to find the action of $H$, since the action of $X_\phi$ is known. The action of $H$ follows from the difference equation \eqref{diff.eq.MPpol} for the Meixner-Pollaczek polynomials.
\begin{prop} \label{prop:rea+-H}
The operator $\Te^+$ intertwines the actions of the generators $H$, $B$, $C$ in the positive discrete series, with the following difference operators:
\[
\begin{split}
\Te^+\pi^+_k(H) &= \Bigg[M\left(\frac{ e^{i\phi}}{i\sin\phi}(k-ix)\right) T_i
+M\left(2\frac{\cos\phi} {\sin\phi}x\right) +M\left(- \frac{e^{-i\phi}}
{i\sin\phi}(k+ix) \right) T_{-i}\Bigg] \Te^+, \\
\Te^+\pi^+_k(B) &= \Bigg[M\left(\frac{e^{2i\phi}}{2i\sin\phi}(k-ix)\right)T_i +
M\left(\frac{x}{\sin\phi}\right) +M\left(-\frac{e^{-2i\phi}}{2i\sin\phi}(k+ix)\right)T_{-i}\Bigg] \Te^+, \\
\Te^+\pi^+_k(C) &= \Bigg[M\left(-\frac{1}{2i\sin\phi}(k-ix)\right)T_i +
M\left(-\frac{x}{\sin\phi}\right) + M\left(\frac{1}{2i\sin\phi}(k+ix)\right)T_{-i}\Bigg] \Te^+,
\end{split}
\]
where $T$ denotes the shift operator: $T_a f(x)=f(x+a)$.
For the negative discrete series, $\Te^-$ intertwines the actions of $H$, $B$, $C$, with the following difference operators:
\[
\begin{split}
\Te^- \pi^-_k(H) &= \Bigg[M\left(-\frac{ e^{i\phi}}{i\sin\phi}(k-ix)\right) T_i
+M\left(-2\frac{\cos\phi} {\sin\phi}x\right) +M\left( \frac{e^{-i\phi}}{i\sin\phi}(k+ix)\right) T_{-i}\Bigg] \Te^-, \\
\Te^- \pi^-_k(B) &=\Bigg[ M\left(-\frac{1}{2i\sin\phi}(k-ix)\right)T_i +
M\left(-\frac{x}{\sin\phi}\right) + M\left(\frac{1}{2i\sin\phi}(k+ix)\right)T_{-i}\Bigg] \Te^-, \\
\Te^- \pi^-_k(C) &= \Bigg[M\left(\frac{e^{2i\phi}}{2i\sin\phi}(k-ix)\right)T_i + M\left(\frac{x}{\sin\phi}\right) +M\left(-\frac{e^{-2i\phi}}{2i\sin\phi}(k+ix)\right)T_{-i}\Bigg] \Te^-.
\end{split}
\]
\end{prop}
\begin{proof}
We find the action of $H$ from the difference
equation \eqref{diff.eq.MPpol} for the Meixner-Pollaczek
polynomials;
\[
\begin{split}
\Te^+\pi^+_k(H)\, e_n &= \Te^+ (2n+2k)\, e_n = (2n+2k)P_n^{(k)}(x;\phi) \\
&=\frac{e^{i\phi}}{i\sin\phi}(k-ix) P_n^{(k)}(x+i;\phi) + 2x\frac{\cos\phi}{\sin\phi}P_n^{(k)}(x;\phi)-\frac{e^{-i\phi}}{i\sin\phi}(k+ix) P_n^{(k)}(x-i;\phi). 
\end{split}
\]
The action of $X_\phi$ is given in Proposition \ref{prop:+-H}: $\Te^+\pi^+_k(X_\phi)\, e_n = 2x \sin\phi\ P_n^{(k)}(x;\phi)$. 
Then Lemma \ref{lem:HXphi} proves the proposition for the positive discrete series. 

We find the action in the negative discrete series in the same way, or we use the Lie-algebra isomorphism
$\thet$, given by
\[
\thet(H) = -H, \quad \thet(B) = C, \quad \thet(C)=B.
\]
Then $\pi^+_k\big(\thet(Y)\big) = \pi^-_k(Y)$ for $Y \in \su(1,1)$.

A straightforward calculation shows that these operators indeed
satisfy the $\su(1,1)$ commutation relations.
\end{proof}
\begin{rem}
To simplify notations we denote $\Te^+ \pi^{+} (\Te^+)^*$ by $\pi^+$, and similarly for $\pi^-$.
\end{rem}
In the same way as in Proposition \ref{prop:rea+-H} it can be shown that $\Te^{\rho,\eps}$ intertwines the actions of $H$, $B$ and $C$ in the principal unitary series with $2 \times 2$ diagonal matrices with difference operators as elements. This is done by finding a difference equation for the Meixner-Pollaczek functions $u_n$ and $u_n^*$, using contiguous relations for $_2F_1$-series. We do not need these realizations here, so we will not work this out.\\

To express $\De(\Om)$ in terms of $H$, $B$ and $C$, the coproduct $\De$ is extended to $\mathcal U\big(\su(1,1) \big)$ as an algebra homomorphism. Then from the definition of the Casimir element \eqref{Casimir} we find
\begin{equation} \label{Copr Om}
\De(\Om) = 1 \tensor \Om + \Om \tensor 1 - \hf H \tensor H - ( C \tensor B + B \tensor C).
\end{equation}
Using this expression and Proposition \ref{prop:rea+-H}, we find the following.
\begin{prop} \label{prop:Om}
In the realizations of Proposition \ref{prop:rea+-H}, we have
\[
\begin{split}
\pi^+_{k_1}& \tensor \pi^-_{k_2}\big( \De(\Om) \big)\Big|_{x_1=x+t,\ x_2=x}=\\
&M\Big(-e^{-2i\phi} \big(k_1+i(t+x)\big)\big(k_2+ix\big)\Big)T_{-i}+ M\Big(k_1(1-k_1) + k_2(1-k_2) - 2(x+t)x\Big) \\ 
+& M\Big(-e^{2i\phi} \big(k_1-i(t+x)\big)\big(k_2-ix\big)\Big)T_{i},
\end{split}
\]
where the shift operator $T$ acts with respect to $x$.
\end{prop}
\begin{proof}
Let $F_1(x)$ and $F_2(x)$ be polynomials in $x$, and let $f(x,t) = F_1(x+t)F_2(x)$, then a large but straightforward computation yields
\[
\begin{split}
\pi^+_{k_1} \tensor \pi^-_{k_2}\big( \De(\Om) \big)f(x,t) =
&[k_1(1-k_1) + k_2(1-k_2) - 2(x+t)x] f(x,t) \\
& -e^{-2i\phi} \big(k_1+i(t+x)\big)\big(k_2+ix\big) f(x-i,t) \\
& -e^{2i\phi} \big(k_1-i(t+x)\big)\big(k_2-ix\big) f(x+i,t).
\end{split}
\]
\end{proof}
\begin{rem}
The action of $\Om$ in the tensor product can also be found from 
\begin{align*}
\pi^+_k(H) &= M(2ix),& \pi^-_k(H) &= M(-2ix),\\
\pi^+_k(B) &= M(k-ix) T_i, &\pi^-_k(B) &= M\left(e^{-2i\phi}(k+ix)\right) T_{-i},\\
\pi^+_k(C) &= M(k+ix) T_{-i}, &\pi^-_k(C) &= M\left(e^{2i\phi}(k-ix)\right) T_{i}.
\end{align*}
These realizations are equivalent to the realizations given in Proposition \ref{prop:rea+-H}.
\end{rem}
In the next section we show that the continuous Hahn functions $\ph_\rho(x;t,k_1,k_2,\phi)$ are eigenfunctions of the difference operator of Proposition \ref{prop:Om}, and we work out the corresponding integral transform.\\

\section{The continuous Hahn integral transform} \label{sec:Hahn}
In this section we study a second order difference operator. This difference operator is obtained from the action of the Casimir operator on hyperbolic basis vectors in the tensor product of a positive and a negative discrete series representation of $\su(1,1)$, see section \ref{ssec:H}. The spectral analysis of this operator, leads to an integral transform pair with a certain type of $_3F_2$-series as a kernel. We call these $_3F_2$-series continuous Hahn functions, because of their similarity to continuous (dual) Hahn polynomials. The method we use is based on asymptotics, and is essentially the same method as used by G\"otze \cite{Go} and Braaksma and Meulenbeld \cite{BM} for the Jacobi function transform by approximating with the Fourier transform.

\subsection{The difference operator $\boldsymbol{\La}$ and the Wronskian}
For $k_1$, $k_2 >0$, $t \in \R$ and $0<\phi<\pi$ the weight function $w(x)$ is defined by
\begin{equation} \label{def:w}
w(x) = \frac{1}{2\pi}e^{(2x+t)(2\phi-\pi)}\Ga(k_1+it+ix) \Ga(k_1-it-ix) \Ga(k_2+ix) \Ga(k_2-ix).
\end{equation}
The difference operator $\La$ is defined by
\[
\La\ :g(x) \mapsto \al_+(x) g(x+i) + \be(x) g(x) + \al_-(x) g(x-i),
\]
where
\[
\begin{split}
\al_\pm(x) &= - e^{\pm 2i\phi}(k_2\mp ix)( k_1\mp i(t+x) ),\\
\be(x) &= k_1(1-k_1) +k_2(1-k_2) -2(t+x)x.
\end{split}
\]
Initially $\La$ is defined for those $g(x)\in L^2(\R, w(x)dx)$ that have an analytic continuation to the strip
\[
\mathcal{S}_\eps =\left\{ z \in \C \ : \  |\Im(z)|<1+\eps,\ \eps>0 \right\}.
\]
The difference operator $\La$ corresponds to the action of the Casimir operator in the tensor product of a positive and a negative discrete series representation on hyperbolic basis vectors, see Proposition \ref{prop:Om}.
\begin{rem} \label{rem:5.1}
Observe that the functions $e^{-x(2\phi-\pi)} p_n(x;k_1+it, k_2, k_1-it, k_2)$, where $p_n$ is a continuous Hahn polynomial in the notation of \cite{KS}, form an orthogonal basis for $L^2(\R, w(x)dx)$, but they are not eigenfunctions of $\La$. Also, the functions $e^{-2\phi x}p_n(x;k_1+it, k_2, k_1-it, k_2)$ are eigenfunctions of $\La$, but they are not elements of $L^2(\R, w(x)dx)$.
\end{rem}

We define 
\[
\inprod{f}{g}_{M,N} = \int_{-M}^N f(x) g^*(x) w(x) dx.
\]
Then $\lim_{N,M \rightarrow \infty} \inprod{f}{g}_{M,N}$
is the inner product on $L^2(\R, w(x)dx)$, so $\inprod{f}{g}_{M,N}$ is a truncated inner product.

\begin{Def} \label{def:Wronsk}
For functions $f$ and $g$ analytic in $\calS_\eps$, the Wronskian $[f,g]$ is defined by
\[
[f,g](y)=\int_y^{y+i} \left\{ f(x) g^*(x-i)- f(x-i) g^*(x) \right\}
\al_-(x) w(x) dx.
\]
\end{Def}
\begin{prop} \label{prop: int=wrons}
Let $f$ and $g$ be analytic in $\calS_\eps$, then
\[
\inprod{\La f}{g}_{N,M} - \inprod{f}{\La g}_{N,M}
= [f,g](N) - [f,g](-M).
\]
\end{prop}
\begin{proof}
Observe that $\al_+^*=\al^-$, $\al_\pm(x \mp i) w(x \mp i) = \al_\pm^*(x) w(x)$ and
$\be^*(x)=\be(x)$. Furthermore for $x\in \R$ we have $g^*(x\pm i)
= \overline { g(x\mp i)}$. This gives for the Wronskian
\[
\begin{split}
[f,g](y) =&  \int_y^{y+i} \left\{ f(x) g^*(x-i)- f(x-i) g^*(x)
\right\} \al_-(x) w(x) dx\\
=& \int_y^{y+i}  f(x) g^*(x-i)\al_+^*(x) w(x) dx - \int_y^{y+i} f(x-i)g^*(x) \al_-(x)w(x) dx \\
=&\int_{y}^{y+i}  f(x) g^*(x-i)\al_+(x-i)w(x-i) dx -\int_{y-i}^y
 f(x) g^*(x+i)\al_-(x+i)w(x+i) dx.
\end{split}
\]
For $f$ and $g$ analytic in $\calS_\eps$, we find
\[
\begin{split}
\int_{-M}^N& (\La f)(x)g^*(x) w(x) dx = \int_{-M}^N
\left[ \al_+(x) f(x+i) + \be(x) f(x) + \al_-(x) f(x-i) \right] g^*(x) w(x) dx \\
=& \int_{-M+i}^{N+i} \al_+(x-i) f(x) g^*(x-i) w(x-i)dx +
\int_{-M}^N \be(x)f(x) g^*(x) w(x) dx \\ &+ \int_{-M-i}^{N-i} \al_-(x+i) f(x) g^*(x+i) w(x+i) dx .
\end{split}
\]
We choose a different path of integration
\[
\int_{-M-i}^{N-i}= \int_{-M-i}^{-M} \ + \int_{-M}^N\ + \int_{N}^{N-i}
\] 
and similarly for $\int_{-M+i}^{N+i}$. Then we find
\[
\begin{split}
\int_{-M}^N& (\La f)(x)g^*(x) w(x) dx = \\
&\int_{-M-i}^{-M} \al_-(x+i) f(x) g^*(x+i)w(x+i) dx -  \int_{N-i}^N \al_-(x+i) f(x) g^*(x+i)w(x+i) dx \\
+& \int_{N}^{N+i} \al_+(x-i) f(x) g^*(x-i)w(x-i) dx -  \int_{-M}^{-M+i} \al_+(x-i) f(x) g^*(x-i)w(x-i) dx\\
+&\int_{-M}^N f(x) \left[ \al_-(x+i) g^*(x+i)w(x+i) +
\be(x) g^*(x)w(x) + \al_+(x-i) g^*(x-i)w(x-i) \right]  dx. 
\end{split}
\]
The last integral $\int_{-M}^N$ equals
\[
\int_{-M}^N f(x)\left[ \al_-^*(x) g^*(x+i) +\be^*(x) g^*(x) +
\al_+^*(x) g^*(x-i) \right] w(x) dx = \int_{-M}^N f(x)  (\La
g)^*(x) w(x)dx,
\]
and in the other four integrals we recognize $[f,g](N) -
[f,g](-M)$.
\end{proof}
Our first goal is to show that $\La$ is a symmetric operator on a domain which will be specified, so we are interested in the limit of the Wronskian $[f,g](y)$ for $y \rightarrow \pm\infty$. The following lemma is useful in determining these limits.
\begin{lem} \label{eq:asymp w}
Let $k_1,k_2>1$, $x\in\R$, and $-1 \leq y \leq 1$, then the weight function $w(x+iy)$ has the following asymptotic behaviour
\[
w(x+iy) =
\begin{cases}
\displaystyle 2\pi e^{2t(\phi-\pi)+4iy(\phi-\pi)} x^{2k_1+2k_2-2} e^{-4(\pi-\phi) x}\left(1+\mathcal
O\big(\frac{1}{x}\big) \right), & x\rightarrow \infty, \\
\displaystyle 2\pi e^{2\phi t+4i\phi y} |x|^{2k_1+2k_2-2} e^{4\phi x}\left(1 +
\mathcal O \big(\frac{1}{x}\big) \right),& x\rightarrow
-\infty.
\end{cases}
\]
\end{lem}
\begin{proof}
From Stirling's asymptotic formula \cite[(8.16)]{Olv} we find for $u>0$ and $v \in \R$
\[
\Ga(u+iv) = \sqrt{2\pi}\, |v|^{u+iv-\hf} e^{-(u+iv)}\, e^{(iu -v)\arctan (v/u)} + \mathcal O\left( \frac{ 1}{v} \right), \qquad v \rightarrow \pm \infty.
\]
We use 
\[
\arctan \frac{v}{u} + \arctan \frac{u}{v} = 
\begin{cases}
\frac{\pi}{2}, & v>0,\\
-\frac{\pi}{2}, & v<0,
\end{cases}
\] 
to find for $v \rightarrow \pm \infty$
\[
(iu-v)\arctan \frac{v}{u} =  (iu-v)\left( \pm\frac{\pi}{2} - \frac{u}{v}\right) + \mathcal O\left(\frac{1}{v} \right).
\]
So we have
\[
\Ga(u+iv) = \sqrt{2\pi}\, |v|^{u+iv-\hf} e^{-\pi|v|/2-iv \pm \pi iu/2} +\mathcal O\left( \frac{ 1}{v} \right), \qquad v \rightarrow \pm \infty.
\]
Applying this formula to the four $\Ga$-functions in \eqref{def:w} gives the asymptotic behaviour of the weight function $w(x+iy)$.
\end{proof}
In general, if for some $\eps>0$ 
\[
\int_\R e^{\eps|x|} d\mu(x) < \infty,
\]
then the moment problem for the measure $d\mu$ is determinate, see e.g.~ \cite{MJ} and references therein. Using this criterion with $0<\eps < \min\{4\phi, 4(\pi-\phi) \}$, we find from Lemma \ref{eq:asymp w} that the moment problem for the measure $w(x)dx$ is determinate. In particular this shows that the polynomials are dense in $L^2(\R, w(x)dx)$.

Let $\D$ be the space of polynomials on $\R$, then $\D$ is a dense subspace of  $L^2(\R, w(x)dx)$.  Since 
\begin{equation} \label{eq:asympt al}
\al_-(x+iy) = e^{-2i\phi} x^2 \Big( 1
+ \mathcal O \big( \frac{1}{x} \big) \Big), \qquad x \rightarrow
\pm \infty,
\end{equation}
it follows from Definition \ref{def:Wronsk} and Lemma \ref{eq:asymp w} that $\lim_{N \rightarrow \pm \infty}[f,g](N)=0$ for $f$, $g$ polynomials. Hence by Proposition \ref{prop: int=wrons} we find
\begin{prop}
The operator $(\La, \D)$ is a densely defined symmetric operator on $L^2(\R, w(x)dx)$.
\end{prop}
\begin{rem}
The operator $(\La, \D)$ is also densely defined and symmetric on the space spanned by $e^{-x(2\phi-\pi)}p_n(x)$, where $p_n$ is a polynomial, cf.~ Remark \ref{rem:5.1}.
\end{rem}

\subsection{Eigenfunctions of $\boldsymbol{\La}$} We determine eigenfunctions of $\La$, using contiguous relations for $_3F_2$-functions. First note that for a monic polynomial of degree $n$, $p_n(x) = x^n + \ldots\ $, we have
\[
(\La p_n)(x) = \left[\al_+(x)+\be(x) + \al_-(x) \right] x^n + \text{lower order terms}.
\]
Since $\al_+(x)+\be(x) + \al_-(x)$ is a polynomial of degree $2$, $\La$ raises the degree of a polynomial by $2$. Therefore $\La$ cannot have polynomial eigenfunctions. 

Let $p(x)$ be the $i$-periodic function
\begin{equation} \label{def: periodic}
p(x)= \frac{1}{\pi}e^{\pi x}\sin(\pi(k_1-it-ix)),
\end{equation}
and let $\ph_\rho(x) = \ph_\rho(x;t,k_1,k_2,\phi)$ and
$\Phi_\rho(x)=\Phi_\rho(x;t,k_1,k_2,\phi)$ denote the functions
\begin{align}
\ph_\rho(x) =& e^{-2\phi x}p(x) \frac{\Ga(1-k_1+it+ix) }{
\Ga(k_1+ix+it) } \F{3}{2}{k_2-ix, k_2-k_1+\hf+i\rho,
k_2-k_1+\hf-i\rho} {2k_2, k_2-k_1+it+1} {1}, \label{def:ph} \\
\Phi_\rho(x) =&e^{-2\phi x}p(x)\frac{\Ga(1-k_1+it+ix) \Ga(1-k_2+ix)}
{\Ga(k_1+ix+it) \Ga(\frac{3}{2}-k_1+ix+i\rho)} \nonumber\\
&\times \F{3}{2}{k_2-k_1+\hf+i\rho, \frac{3}{2}-k_1-k_2+i\rho,
\hf-it+i\rho} {1+2i\rho,
\frac{3}{2}-k_1+ix+i\rho}{1}.\label{def:Phi}
\end{align}
Both $_3F_2$-series are absolutely convergent for $\Re(k_1+it+ix)>0$. Note that the expression for $\ph_\rho(x)$ is the same as \eqref{cont.H.f} after applying Euler's reflection formula.
\begin{prop} \label{prop:eigenf}
For $k_1>1$, the functions $\ph_\rho(x)$ and 
$\Phi_\rho(x)$ are eigenfunctions of $\La$ for eigenvalue
$(\rho^2+\frac{1}{4})$.
\end{prop}
\begin{rem} \label{rem: periodic f1}
(i) Observe that due to the $i$-periodic function $p(x)$, $\ph_\rho(x)$ is an entire function. Denote $\tilde \ph_\rho(x) = \ph_\rho(x)/p(x)$, then $p(x)$ cancels the poles of $\tilde\ph_\rho(x)$. There are more choices of $i$-periodic functions that cancel the poles of $\tilde \ph_\rho(x)$. One of the reasons for this particular choice, is that it appears in the Lie-algebraic interpretation of the function $\ph_\rho(x)$, see Theorem \ref{thm:CG-decomp H}. We come back to the choice of the $i$-periodic function in Remark \ref{rem:period}.

(ii) Observe that $\ph_\rho^*(x)$ is obtained from $\ph_\rho(x)$ by the substitutions $(x,t,\phi) \mapsto (-x, -t, \pi-\phi)$. Since the difference operator $\La$ and the weight function $w(x)$ are invariant under these substitutions, it follows from Proposition \ref{prop:eigenf} that, for $k_1>1$, $\ph_\rho^*(x)$ is also an eigenfunction of $\La$ for eigenvalue $(\rho^2+\frac{1}{4})$. A similar argument shows that $\Phi_\rho^*(x)$ is an eigenfunction of $\La$ for eigenvalue $(\rho^2+\frac{1}{4})$. Obviously $\Phi_{-\rho}(x)$ and $\Phi_{-\rho}^*(x)$ are also eigenfunctions of $\La$ for eigenvalue $(\rho^2+\frac{1}{4})$.
\end{rem}
\begin{proof}
Combining the contiguous relations
\cite[(3.7.9),(3.7.10),(3.7.13)]{AAR}, gives
\[
\begin{split}
(d-a)(e-a) \left[ F(a-)-F \right] + a(a+b+c-d-e+1) \left[F(a+) - F\right] = -bc F,\\
F = \F{3}{2}{a,b,c}{d,e}{1}, \quad F(a\pm) = \F{3}{2}{a\pm 1,b,c}{d,e}{1}.
\end{split}
\]
From this relation we find that
\[
e^{-2\phi x}\frac{\Ga(1-k_1+it+ix) } { \Ga(k_1+ix+it) }
\F{3}{2}{k_2-ix, k_2-k_1+\hf+i\rho, k_2-k_1+\hf-i\rho} { 2k_2,
k_2-k_1+it+1}{1}
\]
is an eigenfunction of $\La $ for eigenvalue $\rho^2+\frac{1}{4}$.
And then, since $p(x)$ is $i$-periodic, $\ph_\rho(x)$ is also an
eigenfunction for eigenvalue $\rho^2+\frac{1}{4}$. Since $\ph_\rho(x)$ must be
analytic in $\calS_\eps$, the condition $k_1>1$ is needed for
absolute convergence of the $_3F_2$-series at the point $x+i$.

Denote
\[
F_{\pm}(a \mp) = \F{3}{2}{a, b \pm 1, c \pm 1}{ d \pm 1, e \pm
1}{1}.
\]
From contiguous relations \cite[(3.7.9), (3.7.10), (3.7.13)]{AAR}
we find
\[
\begin{split}
\Big[ (d-1)(e-1) F_-(a+) - (b-1)(c-1) F \Big] - (a-1)(d+e-a-b-c) F
\\+ \Big[ \frac{ bc(d-a)(e-a)}{de} F_+(a-) - (d-a)(e-a) F \Big]
=0.
\end{split}
\]
From this we see that
\begin{multline} \label{eq:Phi1}
\Phi_\rho(x) = e^{-2\phi x}p(x)\frac{ \Ga(1-k_2+ix) \Ga(1+2i\rho)
\Ga(1- k_1+ix+it) }
{ \Ga(\hf+it+i\rho)\Ga(k_1 +\hf+i\rho +ix) \Ga(\frac{3}{2}-k_1+ix+i\rho) } \\
\times \F{3}{2} { \hf-it+i\rho, 1-k_2+ix, k_2+ix
}{k_1+\hf+ix+i\rho, \frac{3}{2} -k_1 +i\rho +ix }{1},
\end{multline}
is another eigenfunction of $\La$ with eigenvalue
$\rho^2+\frac{1}{4}$. The $_3F_2$ series converges absolutely if
$\Re(\hf+i\rho+it)>0$, so absolute convergence does not depend on
$k_1$. Using \cite[Cor.3.3.5]{AAR} this can be written as
\eqref{def:Phi}, and for this expression the condition $k_1>1$ is
needed. 
\end{proof}

The function $\ph_\rho(x)$ can be expanded in terms of $\Phi_\rho(x)$ and $\Phi_{-\rho}(x)$.
\begin{prop} \label{prop: c-expan}
\[
\ph_\rho(x) = c(\rho)\Phi_\rho(x) + c(-\rho)\Phi_{-\rho}(x),
\]
where
\[
c(\rho) =  \frac{ \Ga(2k_2) \Ga(k_2-k_1+it+1) \Ga(-2i\rho) }{ \Ga(k_1+k_2-\hf-i\rho) \Ga(\hf+it-i\rho)  \Ga(k_2-k_1+\hf-i\rho) }.
\]
\end{prop}
\begin{proof}
This follows from \cite[p.15(2)]{Bai}
\[
\begin{split}
\F{3}{2}{a,b,c}{d,e}{1} =& \frac{ \Ga(1-a) \Ga(d) \Ga(e) \Ga(c-b) }{ \Ga(d-b) \Ga(e-b) \Ga(1+b-a) \Ga(c) } \F{3}{2}{b, b-d+1, b-e+1}{1+b-c, 1+b-a}{1}\\
&+ \text{idem}(b;c).
\end{split}
\]
Here $\text{idem}(b;c)$ after an expression means that the expression is repeated with $b$ and $c$ interchanged. 
\end{proof}
In the next subsection we consider the Wronskians $[\ph_\rho, \ph_\si](y)$ and $[\ph_\rho^*, \ph_\si](y)$ for $y \rightarrow \pm \infty$, so we need the asymptotic behaviour of $\ph_\rho$. We find this from Proposition \ref{prop: c-expan} and the asymptotic behaviour of $\Phi_\rho$.

\begin{lem} \label{lem: assympt}
Let $k_1>1$, $\rho \in \C$ and $-1 \leq y \leq 1$. For $x \rightarrow \pm \infty$
\[
\begin{split}
\Phi_\rho(x+iy) &= p(x+iy) e^{-2\phi(x+iy)} (ix)^{\hf-k_1-k_2-i\rho}
 \left\{1+\frac{1}{2ix}\big(A(\rho)+yB(\rho)\big) +\mathcal
O\big(\frac{1}{x^2} \big) \right\},\\
\Phi_\rho^*(x+iy) &= p^*(x+iy) e^{-2\phi (x+iy)}
(-ix)^{\hf-k_1-k_2+i\overline{\rho}} 
 \left\{1+\frac{1}{-2ix}\big(\overline{A(\rho)}- y\overline{B(\rho)}\big)
+\mathcal O\big(\frac{1}{x^2} \big) \right\},
\end{split}
\]
where
\[
\begin{split}
A(\rho) =& 2it(1-2k_1) + (k_1-k_2-\hf-i\rho)(\frac{3}{2}-k_1-k_2+i\rho) \\
&+ \frac{(k_2-k_1+\hf+i\rho)(\frac{3}{2}-k_1-k_2+i\rho)(\hf-it+i\rho)} {\hf+i\rho},\\
B(\rho) =& 2k_1+2k_2-1+2i\rho.
\end{split}
\]
\end{lem}
\begin{proof}
This follows from \eqref{def:Phi} and the asymptotic formula for the ratio of two $\Ga$-functions \cite[\S4.5]{Olv}
\begin{equation}\label{assymp Ga}
\frac{ \Ga (a+z)}{ \Ga(b+z) } = z^{a-b}\left( 1+
\frac{1}{2z}(a-b)(a+b-1) + \mathcal O \big(\frac{1}{z^2}\big) \right), \quad |z|\rightarrow \infty, \quad |\arg(z)|<\pi.
\end{equation}
The first part of the expression for $A(\rho)$ is obtained from
\[
(a_1-b_1)(a_1+b_1-1) + (a_2-b_2)(a_2+b_2-1),
\]
the second part comes from the second term in the hypergeometric series. 
\end{proof}
The asymptotic behaviour of the $i$-periodic function $p(x)$ is also needed;
\[
p(x+iy)=\ \frac{ e^{i \pi (k_1+2y)}}{2\pi i} e^{\pi(2x+t)} - \frac{e^{-i \pi k_1}}{2\pi i} e^{-\pi t} 
= 
\begin{cases}
\displaystyle \frac{ e^{i \pi (k_1+2y)}}{2\pi i} e^{\pi(2x+t)} + \mathcal O(1), & x \rightarrow \infty,\\ \\
\displaystyle - \frac{e^{-i \pi k_1}}{2\pi i} e^{-\pi t} + \mathcal O (e^{2\pi x}),& x \rightarrow -\infty,
\end{cases}
\]
\[
p^*(x+iy)=\ \frac{ e^{i \pi k_1}}{2\pi i} e^{-\pi t} - \frac{e^{-i \pi (k_1-2y)}}{2\pi i} e^{\pi(2x+t)} 
= 
\begin{cases}
\displaystyle -\frac{ e^{-i \pi (k_1-2y)}}{2\pi i} e^{\pi(2x+t)}  + \mathcal O(1), & x \rightarrow \infty,\\ \\
\displaystyle  \frac{e^{i \pi k_1}}{2\pi i} e^{-\pi t} + \mathcal O (e^{2\pi x}),& x \rightarrow -\infty.
\end{cases}
\]

\subsection{Continuous spectrum} \label{ssect: continuous}
We determine the spectrum of the difference operator $\La$. In this subsection we consider the case where the spectrum only consists of a continuous part. 

Since
\[
|(ix)^{\hf-k_1-k_2-i\rho}|^2 =
\begin{cases}
|x|^{1-2k_1-2k_2+2\Im(\rho)}e^{\pi \Re(\rho)}, & x>0,\\
|x|^{1-2k_1-2k_2+2\Im(\rho)}e^{-\pi \Re(\rho)}, & x<0,
\end{cases}
\]
we find from Proposition \ref{lem: assympt}
\begin{equation} \label{Phi w}
|\Phi_\rho(x)|^2 w(x) = \mathcal O (|x|^{2\Im(\rho)-1}), \qquad x \rightarrow \pm \infty.
\end{equation}
So $\Phi_\rho(x)$ is an element of $L^2(\R, w(x)dx)$ for $\Im(\rho)<0$. This shows that it is possible to give eigenfunctions of $\La$ for complex eigenvalues. 
We only consider eigenfunctions which are even in $\rho$, and in that case all eigenvalues of $\La$ are real.   

First we consider the continuous spectrum of $\La$. We show that $[\frac{1}{4}, \infty)$ is contained in the continuous spectrum. 
Assume that $\rho$ is real and that the $c$-function in Proposition \ref{prop: c-expan} does not have zeros, or, equivalently, assume that $k_1+k_2 \geq \hf$ and $k_2-k_1 > -\hf$. Since we only consider even functions in $\rho$, we may assume $\rho \geq 0$. We use Proposition \ref{prop: int=wrons} to calculate the truncated inner product of two eigenfunctions. This gives for $k_1>1$ and $\rho \neq \si$
\begin{equation} \label{inp-Wr}
\inprod{\ph_\rho}{ \ph_\si}_{M,N} =\frac{ [\ph_\rho, \ph_\si](N) - [\ph_\rho, \ph_\si](-M)}{\rho^2-\si^2}.
\end{equation}
Multiplying both sides with an arbitrary function $f(\rho)$ and 
integrating over $\rho$ from $0$ to $\infty$, gives
\begin{equation} \label{eq:Wr}
\int_0^\infty f(\rho) \inprod{\ph_\rho}{ \ph_\si}_{M,N}\ d\rho = \int_0^\infty f(\rho) \frac{ [\ph_\rho, \ph_\si](N)
- [\ph_\rho, \ph_\si](-M) }{\rho^2 - \si^2} d\rho.
\end{equation}
The function $f(\rho)$ must satisfy some conditions, which we
shall determine later on. We take limits $N,M \rightarrow \infty$
on both sides. To determine the limits of the Wronskians, the following lemma is used. 

\begin{lem} \label{lem:2 app}
Let $k_1,k_2>1$ and $\rho,\si \geq 0$. For for $x \rightarrow \pm\infty$
\[
[\Phi_\rho, \Phi_\si](x) = \pm D^{\pm}(\rho,\si)(\rho+\si)
 |x|^{i(\si-\rho)}\Big( 1 + \mathcal O(\frac{1}{x}) \Big),
\]
where
\[
D^{\pm}(\rho,\si)= \frac{i}{2\pi}e^{ t(2\phi-\pi)} e^{\pm\hf\pi(\rho+\si+2t)} .
\]
\end{lem}
\begin{proof}
From Lemma \ref{lem: assympt} we find for $0 \leq y \leq 1$ and $x \rightarrow \infty$
\[
\begin{split}
\Phi_\rho(x+iy) &\Phi_\si^*(x+iy-i)-\Phi_\rho(x+iy-i) \Phi_\si^*(x+iy) = \\
&\, e^{2i\phi-4i\phi y+4i\pi y -4\phi x} |x|^{1-2k_1-2k_2}  |p(x)|^2 (ix)^{-i\rho}(-ix)^{i\si} \frac{B(\rho) -\overline{B(\si)}   }{2ix}\Big( 1 + \mathcal O(\frac{1}{x}) \Big) \\
=&\,e^{2i\phi-4i\phi y+4i\pi y -4\phi x} |x|^{1-2k_1-2k_2}  |p(x)|^2 e^{\pm \hf \pi(\rho+\si) } |x|^{i(\si-\rho)} \frac{\rho+\si}{x}\Big( 1 + \mathcal O(\frac{1}{x}) \Big).
\end{split}
\]
Using the asymptotic behaviour of $p(x)$ and Lemma \ref{eq:asymp w}, we find the asymptotic behaviour of the integrand of the Wronskian. Note that this is independent of $y$. In a similar way we find the same asymptotic behaviour of the integrand for $x\rightarrow -\infty$.  Now the lemma follows from writing the Wronkian as
\[
\int_x^{x+i} f(z) dz = i\int_0^1 f(x+is) ds,
\] 
and applying dominated convergence.
\end{proof}

\begin{prop} \label{prop:3}
Let $\si \geq 0$, and let $f$ be a continuous function satisfying
\[
f(\rho) =
\begin{cases}
\mathcal O( e^{-\pi \rho}\rho^{2k_2-\hf-\eps}), & \rho \rightarrow \infty, \quad \eps>0, \\
\mathcal O(\rho^\de), & \rho \rightarrow 0, \quad \de>0,
\end{cases}
\]
then
\[
\lim_{N\rightarrow \infty}\int_0^\infty f(\rho) \inprod{\ph_\rho}{ \ph_\si}_{N,N}\, d\rho = W_0^{-1}(\si) f(\si), \]
where
\[
W_0(\si) = 
 \frac{ 1}{2\pi} e^{- t(2\phi-\pi)}\, \left| \frac{ \Ga( k_2-k_1+\hf+i\si) \Ga(k_1+k_2-\hf+i\si)\Ga(\hf+it+i\si) \Ga(\hf-it+i\si) }{ \Ga(2k_2) \Ga(k_2-k_1+it+1) \Ga(2i\si) }\right|^2 .
\]
\end{prop}
\begin{proof}
Let $k_1,k_2>1$. We use \eqref{eq:Wr} to calculate $\lim_{N \rightarrow \infty}\int_0^\infty f(\rho) \inprod{\ph_\rho}{ \ph_\si}_{N,N}\ d\rho$. From the $c$-function expansion, see Proposition \ref{prop: c-expan}, we obtain
\[
[\ph_\rho, \ph_\si](x) = \sum_{\epsilon, \xi \in \{-1,1\}} 
c( \epsilon \rho) \overline{c(\xi \si)} [\Phi_{\epsilon\rho}, \Phi_{\xi\si}](x).
\]
Then \eqref{inp-Wr} and Lemma \ref{lem:2 app} give, for $N \rightarrow \infty$,
\[
\begin{split}
\inprod{\ph_\rho}{\ph_\si}_{N,N} = \sum_{\epsilon, \xi \in \{-1,1\}}\frac{N^{i(\xi \si - \epsilon \rho)}}{\epsilon \rho - \xi \si} c( \epsilon \rho) \overline{c(\xi \si)} \Big[D^+(\epsilon \rho, \xi \si) + D^-(\epsilon \rho, \xi \si) \Big]  \Big(1+ \mathcal O\big(\frac{1}{N}\big) \Big).
\end{split}
\]
From \eqref{eq:Wr} we find
\[
\begin{split}
\lim_{N \rightarrow \infty} \int_0^\infty f(\rho) &\inprod{\ph_\rho}{ \ph_\si}_{N,N}\ d\rho =\\
\lim_{N \rightarrow \infty} \int_0^\infty f(\rho)
&\Bigg\{ \psi_1(\rho) \cos\left([\rho+\si] \ln N \right) + \psi_2(\rho) \sin\left([\rho+\si] \ln N \right) \\
&+\psi_3(\rho) \cos\left([\rho-\si] \ln N \right) +
\psi_4(\rho) \frac{\sin\left([\rho-\si] \ln N \right)} {\rho-\si} \Bigg\}
d\rho,
\end{split}
\]
where
\[
\begin{split}
\psi_1(\rho) &= \frac{1} {\rho+\si} \left[ c(\rho) \overline{ c(-\si) } \big(D^+(\rho,-\si) + D^-(\rho,-\si) \big) - c(-\rho) \overline{ c(\si) } \big(D^+(-\rho,\si) + D^-(-\rho,\si) \big) \right],\\
\psi_2(\rho) &= \frac{-i} {\rho+\si} \left[ c(\rho) \overline{ c(-\si) } \big(D^+(\rho,-\si) + D^-(\rho,-\si) \big) + c(-\rho) \overline{ c(\si) } \big(D^+(-\rho,\si) + D^-(-\rho,\si) \big) \right],\\
\psi_3(\rho) &= \frac{1} {\rho-\si} \left[ c(\rho) \overline{ c(\si) } \big(D^+(\rho,\si) + D^-(\rho,\si) \big) - c(-\rho) \overline{ c(-\si) } \big(D^+(-\rho,-\si) + D^-(-\rho,-\si) \big) \right],\\
\psi_4(\rho) &= -i \left[ c(\rho) \overline{ c(\si) } \big(D^+(\rho,\si) + D^-(\rho,\si) \big) + c(-\rho) \overline{ c(-\si) } \big(D^+(-\rho,-\si) + D^-(-\rho,-\si) \big) \right].
\end{split}
\]
Writing out explicitly the terms between square brackets for $\psi_3$ gives, for $\rho=\si$,
\[
e^{t(2\phi-\pi)}\left| \frac{ \Ga(2k_2) \Ga(k_2-k_1+it+1) \Ga(2i\rho)} { \Ga(k_1+k_2-\hf+i\rho) \Ga(k_2-k_1+\hf+i\rho)  } \right|^2\Bigg( \frac{\sin\big( \pi( \hf +it +i\rho) \big) }{ |\Ga(\hf-it+i\rho)|^2} -  \frac{\sin\big( \pi( \hf -it +i\rho) \big)}{ |\Ga(\hf+it+i\rho)|^2} \Bigg).
\]
From Euler's reflection formula it follows that this is equal to zero. So $\psi_3$ has a removable singularity at the point $\rho=\si$.

From the Riemann-Lebesgue lemma, see e.g.~ \cite[\S 9.41]{WW}, it follows 
that for $f\psi_i \in L^1(0,\infty)$, $i=1,2,3$, the terms with $\psi_i$, $i=1,2,3$, vanish. This leaves us with a Dirichlet integral, for which we have the property (see e.g. \cite[\S 9.7]{WW})
\begin{equation} \label{Dir kern}
\lim_{t \rightarrow \infty}\frac{1}{\pi} \int_0^\infty g(x) \frac{
\sin[t(x-y)]}{x-y} dx = g(y) , 
\end{equation}
for a continuous function $g \in L^1(0,\infty)$.
This gives for a continuous function $f$ that satisfies $f\psi_4^\pm \in L^1(0,\infty)$,
\[
\begin{split}
\lim_{N\rightarrow \infty}&\int_0^\infty f(\rho)\inprod{\ph_\rho}{ \ph_\si}_{N,N}\ d\rho = \pi f(\si)\psi_4(\si) \\
=&\ \hf e^{ t(2\phi-\pi)}\, \left\{(e^{-\pi (\si+t)}+e^{\pi (\si+t)}) c(\si) \overline{c(\si)} +(e^{-\pi (\si-t)}+e^{\pi (\si-t)})c(-\si) \overline{c(-\si) } \big)\right\} f(\si)\\
=&\  e^{ t(2\phi-\pi)}\, \left| \frac{ \Ga(2k_2)
\Ga(k_2-k_1+it+1) \Ga(2i\si) } { \Ga( k_2-k_1+\hf+i\si)
\Ga(k_1+k_2-\hf+i\si) } \right|^2 \\
& \times \Bigg( \frac{\sin\big( \pi( \hf +it +i\si) \big) }{ |\Ga(\hf-it+i\si)|^2} +  \frac{\sin\big( \pi( \hf -it +i\si) \big)}{ |\Ga(\hf+it+i\si)|^2} \Bigg)f(\si)\\
=&\ 2\pi e^{ t(2\phi-\pi)}\, \left| \frac{ \Ga(2k_2) \Ga(k_2-k_1+it+1) \Ga(2i\si) } { \Ga( k_2-k_1+\hf+i\si) \Ga(k_1+k_2-\hf+i\si)\Ga(\hf+it+i\si) \Ga(\hf-it+i\si) } \right|^2 f(\si)\\
=& W_0^{-1}(\si) f(\si).
\end{split}
\]
In the last step Euler's reflection formula is used. The conditions $k_1,k_2>1$ can be removed by analytic continuation. Since
\[
W_0^{-1}(\si) =
\begin{cases}
\mathcal O(\si^{1-4k_2} e^{2\pi \si}), & \si\rightarrow \infty,\\
\mathcal O(1), & \si \rightarrow 0,
\end{cases}
\]
we find that the proposition is valid for the conditions of $f$ as stated.
\end{proof}

Next we consider the truncated inner product $\inprod{\ph_\rho^*}{\ph_\si}_{M,N}$ and the corresponding Wronskians. We need to find the analogues of Lemma \ref{lem:2 app} and Proposition \ref{prop:3}.
\begin{lem} \label{lem:2a}
Let $k_1>1$ and $\rho,\si \geq 0$. For $x \rightarrow \pm\infty$
\[
[\Phi_\rho^*, \Phi_\si](x) = E^{\pm}(\rho,\si)(\rho-\si)
 |x|^{i(\rho+\si)}\Big( 1 + \mathcal O(\frac{1}{x}) \Big),
\]
where
\[
E^{\pm}(\rho,\si)= \frac{1}{2\pi}e^{ t(2\phi-\pi)} e^{\pm \hf\pi i(2k_1-2k_2+i\rho+i\si+2it)} .
\]
\end{lem}
\begin{proof}
The proof is similar to the proof of Lemma \ref{lem:2 app}.
\end{proof}
\begin{prop} \label{prop:3a}
Let $\si \geq 0$, and let $f$ be an even continuous function satisfying
\[
f(\rho) =
\begin{cases}
\mathcal O( e^{-\pi \rho}\rho^{2k_2-\hf-\eps}), & \rho \rightarrow \infty, \quad \eps>0, \\
\mathcal O(\rho^\de), & \rho \rightarrow 0, \quad \de>0,
\end{cases}
\]
then
\[
\lim_{N\rightarrow \infty}\int_0^\infty f(\rho) \inprod{\ph_\rho^*}{ \ph_\si}_{N,N}\, d\rho = W_1^{-1}(\si) f(\si), 
\]
where
\[
\begin{split}
W_1(\si) = &
 \frac{1}{2\pi} e^{- t(2\phi-\pi)}\, \frac{ \Ga(k_1-k_2+it) \Ga(\hf-it-i\si) \Ga(\hf-it+i\si)}{ \Ga(k_2-k_1-it+1) \Ga(2k_2)^2 } \\
&\times \left| \frac{ \Ga( k_2-k_1+\hf+i\si) \Ga(k_1+k_2-\hf+i\si) }{  \Ga(2i\si) }\right|^2 .
\end{split}
\]
\end{prop}
\begin{proof}
The proof runs along the same lines as the proof of Proposition \ref{prop:3}, therefore we leave out the details.

As in the proof of Proposition \ref{prop:3} we find from Lemma \ref{lem:2a}
\[
\begin{split}
\lim_{N\rightarrow \infty} \int_0^\infty f(\rho)& \inprod{\ph_\rho^*}{\ph_\si}_{N,N}\, d\rho =\\
\lim_{N \rightarrow \infty} \int_0^\infty f(\rho)
&\Bigg\{ \psi_1(\rho) \cos\left([\rho+\si] \ln N \right) + \psi_2(\rho) \sin\left([\rho+\si] \ln N \right) \\
&+\psi_3(\rho) \cos\left([\rho-\si] \ln N \right) +
\psi_4(\rho) \frac{\sin\left([\rho-\si] \ln N \right)} {\rho-\si} \Bigg\}
d\rho,
\end{split}
\]
where
\[
\begin{split}
\psi_1(\rho) &= \frac{1} {\rho+\si} \left[  \overline{c(\rho)} \overline{ c(\si) } \big(E^+(\rho,\si) - E^-(\rho,\si) \big) -  \overline{c(-\rho)} \overline{ c(-\si) } \big(E^+(-\rho,-\si) - E^-(-\rho,-\si) \big) \right],\\
\psi_2(\rho) &= \frac{i} {\rho+\si} \left[  \overline{c(\rho)} \overline{ c(\si) } \big(E^+(\rho,\si) - E^-(\rho,\si) \big) +  \overline{c(-\rho)} \overline{ c(-\si) } \big(E^+(-\rho,-\si) - E^-(-\rho,-\si) \big) \right],\\
\psi_3(\rho) &= \frac{-1} {\rho-\si} \left[ \overline{ c(-\rho)} \overline{ c(\si) } \big(E^+(-\rho,\si) - E^-(-\rho,\si) \big) -  \overline{c(\rho)} \overline{ c(-\si) } \big(E^+(\rho,-\si) - E^-(\rho,-\si) \big) \right],\\
\psi_4(\rho) &= -i \left[ \overline{ c(-\rho) }\overline{ c(\si) } \big(E^+(-\rho,\si) - E^-(-\rho,\si) \big) +  \overline{c(\rho)} \overline{ c(-\si) } \big(E^+(\rho,-\si) - E^-(\rho,-\si) \big) \right].
\end{split}
\]
For $\rho=\si$ the term between the square brackets for $\psi_3$ is equal to zero, so $\psi_3$ has a removable singularity at the point $\rho=\si$. So, for $f\psi_i \in L^1(0,\infty)$, $i=1,2,3$, the terms with $\psi_{i}$, $i=1,2,3$, vanish by the Riemann-Lebesgue lemma. This leaves us with a Dirichlet integral. Then, after applying \eqref{Dir kern}, we find for $f\psi_4 \in L^1(0,\infty)$
\[
\lim_{N\rightarrow \infty} \int_0^\infty f(\rho) \inprod{\ph_\rho^*}{\ph_\si}_{N,N}\, d\rho =\pi \psi_4(\si) f(\si).
\]
Writing out $\psi_4$ explicitly gives the result.
\end{proof}
\begin{rem} \label{rem:period}
In Remark \ref{rem: periodic f1} we observed that the $i$-periodic function $p(x)$ cancels the poles of $\tilde \ph_\rho(x)$. Other obvious choices with the same property would be $e^{2k\pi x}p(x)$, for $k \in \Z$. However for $k \neq 0$, the method we used here to find an integral transform pair would fail, since the method depends on the use of the Riemann-Lebesgue lemma and the Dirichlet kernel, which can no longer be used in case $k \neq 0$. This can e.g.~ be seen from Lemma \ref{lem:2 app}, where the terms in front of $|x|^{i(\si-\rho)}$ would contain a factor $e^{2k\pi x}$. So this gives a heuristic argument for the choice \eqref{def: periodic} of the $i$-periodic function.
\end{rem}

Let $f$ be a continuous function satisfying, 
\begin{equation} \label{cond f}
f(\rho) =
\begin{cases}
\mathcal O(\rho^{-2k_2-\eps}e^{\pi \rho} ), & \rho \rightarrow \infty, \quad \eps>0,\\
\mathcal O(\rho^\de), & \rho \rightarrow 0, \quad \de>0,
\end{cases}
\end{equation}
and let $\mathbf{f}$ be the vector
\[
\mathbf{f}(\rho) = \vect{\overline{f(\rho)}}{f(\rho)}.
\]
We define, for $x \in \R$, an operator $\tilde{\mathcal F}$ by
\begin{equation} \label{def:F}
(\tilde{\mathcal F} \mathbf{f})(x) = \int_0^\infty \Big( \ph_\rho(x) \overline{f(\rho)}+\ph_\rho^*(x)f(\rho)\Big) W_0(\rho)  d\rho.
\end{equation}
To verify that this is a well-defined expression, we determine the behaviour of $\ph_\rho(x)$ and $W_0(\rho)$ for $\rho \rightarrow \infty$ and $\rho \downarrow 0$. From Thomae's transformation \cite[Cor.3.3.6]{AAR} we find
\[
\ph_\rho(x)=
\begin{cases}
\mathcal O (\rho^{-2k_2} e^{\pi \rho} ), &\rho \rightarrow \infty,\\
\mathcal O(1), & \rho \rightarrow 0.
\end{cases}
\]
And from \eqref{assymp Ga} we obtain
\begin{equation} \label{assym W}
W_0(\rho) =
\begin{cases}
\mathcal O(\rho^{4k_2-1} e^{-2\pi \rho}), & \rho \rightarrow \infty,\\
\mathcal O(1), & \rho \rightarrow 0.
\end{cases}
\end{equation}
Then we see that the integral in \eqref{def:F} converges absolutely for $f$ satisfying \eqref{cond f}.

For a continuous function $g$ satisfying 
\begin{equation} \label{cond g}
g(x) = 
\begin{cases}
\mathcal O\big(|x|^{\hf-k_1-k_2-\eps}e^{2(\pi-\phi)x}\big), & x \rightarrow \infty , \quad \eps>0,\\
\mathcal O\big(|x|^{\hf-k_1-k_2-\de}e^{-2\phi x}\big), & x \rightarrow -\infty , \quad \de>0,
\end{cases}
\end{equation}
we define, for $\rho \geq 0$, an operator $\mathcal G$ by
\begin{equation} \label{def G}
(\mathcal G g)(\rho) = \int_\R g(x) \vect{\ph_\rho^*(x)}{\ph_\rho(x)} w(x) dx .
\end{equation}
From the asymptotic behaviour of $\ph_\rho^*(x)$ and $w(x)$ for $x \rightarrow \pm \infty$, see Lemmas \ref{eq:asymp w} and \ref{lem: assympt}, it follows that the integral in \eqref{def G} converges absolutely.

\begin{prop} \label{prop:F inv}
If $g=\tilde{\mathcal F} \mathbf{f}$, and $g$ satisfies the conditions \eqref{cond g}, then 
\[
(\mathcal G g)(\rho)= \begin{pmatrix} 1 & W_0(\rho)/W_1(\rho) \\ W_0(\rho)/\overline{W_1(\rho)} & 1 
\end{pmatrix}
\vect{\overline{f(\rho)}}{f(\rho)} .
\]
\end{prop}
\begin{proof}
For a function $g$ satisfying \eqref{cond g} we define operators $\mathcal G_1$ and $\mathcal G_2$ by
\begin{equation} \label{G1 G2}
\begin{split}
(\mathcal G_1 g)(\rho) &= \int_\R g(x) \ph_\rho^*(x) w(x) dx, \\
(\mathcal G_2 g)(\rho) &= \int_\R g(x) \ph_\rho(x) w(x) dx,
\end{split}
\end{equation}
then we have
\[
(\mathcal G g)(\rho) = \vect{ (\mathcal G_1g)(\rho) }{ (\mathcal G_2g)(\rho) }.
\]
If $g(x)=(\tilde{\mathcal F} f)(x)$ satisfies the conditions \eqref{cond g}, then the integral
\[
(\mathcal G_1g)(\si)= \int_\R (\tilde{\mathcal F} f)(x) \ph_\si^*(x) w(x) dx
\]
converges absolutely. So from \eqref{def:F} we obtain
\[
(\mathcal G_1g)(\si)= \lim_{N \rightarrow \infty} \int_{-N}^N \ph_\si^*(x) \left\{\int_0^\infty \ph_\rho(x) \overline{f(\rho)}  W_0(\rho) d\rho +\int_0^\infty \ph_\rho^* (x) f(\rho) W_0(\rho) d\rho\right\}  w(x)dx,
\]
and interchanging integration gives
\[ 
(\mathcal G_1g)(\si) = \lim_{N \rightarrow \infty} \int_0^\infty \overline{f(\rho)} W_0(\rho) \inprod{\ph_\rho}{\ph_\si}_{N,N}\, d\rho + \lim_{N \rightarrow \infty} \int_0^\infty f(\rho) W_0(\rho) \inprod{\ph_\rho^*}{\ph_\si}_{N,N}\, d\rho.
\]
From \eqref{cond f} and \eqref{assym W} it follows that the functions $\overline{f(\rho)}W_0(\rho)$, $f(\rho) W_0(\rho)$ satisfy the conditions for Propositions \ref{prop:3}, \ref{prop:3a} respectively. Applying the propositions gives $(\mathcal G_1 g)(\si) = \overline{f(\si)} + f(\si) W_0(\si)/ W_1(\si)$. In the same way $(\mathcal G_2 g)(\si)$ can be calculated. So we find
$(\mathcal G_2 g)(\si)=f(\si) + \overline{f(\si)} W_0(\si)/\overline{ W_1(\si)}$. 
\end{proof}

We define an operator $\mathcal F$ by
\begin{equation} \label{def F}
\begin{split}
(\mathcal F \mathbf{f})(x) &= \left(\tilde{\mathcal F} \begin{pmatrix} 1 & W_0/W_1 \\ W_0/\overline{W_1} & 1 
\end{pmatrix}^{-1}
\mathbf{f}\right)(x) \\
&= \int_0^\infty  \vect{\ph_\rho^*(x)}{\ph_\rho(x)}^* 
\begin{pmatrix} 
1 & - W_0(\rho)/W_1(\rho) \\ - W_0(\rho)/\overline{W_1(\rho)} &1
\end{pmatrix}
\vect{\overline{f(\rho)}}{f(\rho)}\frac{W_0(\rho)}{1- \left|\frac{W_0(\rho)}{W_1(\rho)} \right|^2} d\rho.
\end{split}
\end{equation}
From Proposition \ref{prop:F inv} we find $\big(\mathcal G (\mathcal F \mathbf {f})\big)(\rho) = \mathbf{f}(\rho)$.

\begin{rem} \label{rem:psi}
From Euler's reflection formula and the identity $2\sin x \sin y = \cos(x-y)-\cos(x+y)$, we find
\[
1-\left| \frac{W_0(\rho)}{W_1(\rho)} \right|^2 = \left|\frac{ \Ga( \hf +it +i\rho) \Ga( \hf-it +i\rho) }{\Ga(k_1-k_2+\hf+i\rho) \Ga(k_2-k_1+\hf+i\rho) } \right|^2.
\]
We define 
\[
\begin{split}
\psi_\rho(x) = & \frac{e^{-x(2\phi-\pi)}\, \Ga(\hf-it-i\rho)  \Ga(\hf-it+i\rho)\Ga(2k_2) \Ga(k_2-k_1+it+1) }{ \Ga(k_1+ix+it) \Ga(k_2-ix) \Ga(k_1+k_2-it) \Ga(k_1-k_2-it+1)}\\ 
&\times \F{3}{2}{ \hf-it+i\rho, \hf-it-i\rho, k_1-ix-it}{ k_1-k_2-it+1, k_1+k_2-it}{1},
\end{split}
\]
then from \cite[\S3.8(1)]{Bai} we obtain
\[
\ph_\rho(x)- \frac{W_0(\rho)}{\overline{W_1(\rho)}}\, \ph_\rho^*(x) = \frac{ \psi_\rho(x)}{|\Ga(k_2-k_1+\hf+i\rho)|^2}.
\]
So the definition of $\mathcal F$ \eqref{def F} is equivalent to
\[
(\mathcal F \mathbf{f})(x) = \int_0^\infty \vect{\psi_\rho^*(x)}{\psi_\rho(x)}^* \vect{ \overline{f(\rho)}}{f(\rho)} W_2(\rho) d\rho,
\]
where 
\[
W_2(\rho) = \frac{e^{-t(2\phi-\pi)}}{2\pi} \left| \frac{ \Ga(k_1-k_2+\hf+i\rho) \Ga(k_2-k_1+\hf+i\rho) \Ga(k_1+k_2-\hf+i\rho) }{ \Ga(2k_2) \Ga(k_2-k_1+it+1) \Ga(2i\rho)} \right|^2 .
\]
Note that $W_2(\rho)d\rho$ is the orthogonality measure for the continuous dual Hahn polynomials. 
\end{rem}

The function $\mathcal F \mathbf{f}$ exists for all functions $\mathbf{f}$ for which the integral \eqref{def F} converges. We want to find a domain on which $\mathcal F$ is injective and isometric. We look for a set $\mathcal S$ of functions for which $\mathcal F \mathcal S$ is a dense subspace of $L^2(\R, w(x)dx)$. Recall that the set of polynomials is a dense subspace of $L^2(\R, w(x)dx)$.

\begin{lem} \label{MP trans}
Let $p_n^{(\la)}(\cdot;\phi)$ denote a Meixner-Pollazcek polynomial as defined by \eqref{def:MP pol}, then
\[
\big(\mathcal G p^{(k_2)}_n(\cdot;\phi) \big) (\rho) = \mathbf{q_n}(\rho) 
= \vect{\overline{q_n(\rho)}}{q_n(\rho)},
\]
with $q_n$ given by
\[
\begin{split}
q_n(\rho) =&e^{t(2\phi - \pi)}  \left( \frac{ -e^{-i\phi} }{ 1-e^{-2i\phi}} \right)^{2k_2+n}  \frac{(-1)^n \Ga(2k_2)\,|(k_2-k_1+\hf+i\rho)_n|^2 }{n!\, (k_2-k_1+it+1)_n } \\
 &\times  \F{2}{1}{n+k_2-k_1+\hf+i\rho, n+k_2-k_1+\hf-i\rho }{ n+k_2-k_1+it+1}{ \frac{1}{ 1- e^{2i\phi} }}. 
\end{split}
\]
\end{lem}
\begin{proof}
Let $\mathcal G_1$ and $\mathcal G_2$ be as defined by \eqref{G1 G2}. Since the Meixner-Pollaczek  polynomial $p_n$ is real, we have $\mathcal G_1 p_n = \overline{ \mathcal G_2 p_n}$. Writing out $\overline{q_n(\rho)}=\left(\mathcal G_1 p^{(k_2)}_n(\cdot;\phi) \right)(\rho)$ explicitly, gives
\[
\begin{split}
\overline{q_n(\rho)}= e^{t(2\phi - \pi)} e^{-in\phi} \frac{ (2k_2)_n }{n!} \frac{1}{2\pi}\int_\R& e^{x(2\phi-\pi)} | \Ga(k_2+ix) |^2 \F{2}{1}{-n, k_2-ix}{2k_2}{1-e^{2i\phi}}\\
 \times& \F{3}{2}{ k_2+ix, k_2-k_1+\hf+i\rho, k_2-k_1+\hf-i\rho}{ 2k_2, k_2-k_1-it+1}{1}  dx \\
=e^{t(2\phi - \pi)} e^{-in\phi} \frac{ (2k_2)_n }{n!} \sum_{m=0}^\infty \sum_{l=0}^n&  \frac{(k_2 - k_1+\hf+i\rho)_m (k_2 - k_1+\hf-i\rho)_m }{m!\, (2k_2)_m (k_2-k_1-it+1)_m} \frac{ (-n)_l }{l!\, (2k_2)_l } (1-e^{2i\phi})^l \\
 \times \frac{1}{2\pi}&\int_\R e^{x(2\phi - \pi)} \Ga(k_2+ix+m) \Ga(k_2-ix+l) dx.
\end{split}
\] 
The inner integral can be evaluated by \cite[(3.3.9)]{PK},
\[
\frac{1}{2\pi i } \int_{k-i\infty}^{k+i\infty} \Ga(s) \Ga(a-s) y^{-s} ds = \frac{ \Ga(a) }{ (1+y)^a }, \qquad 0<\Re(s)<\Re(a), \quad |\arg(y)|< \pi,
\]
with $s=k_2+m+ix$ and $a=2k_2+m+l$. Now the sum over $l$ becomes a terminating $_2F_1$-series, which can be evaluated by the Chu-Vandermonde identity \cite[Cor.2.2.3]{AAR}
\[
\F{2}{1}{-n, 2k_2+m}{2k_2}{1} = 
\begin{cases}
0, & m<n,\\
\displaystyle \frac{(-m)_n }{(2k_2)_n }, & m \geq n.
\end{cases}
\]
Then $\overline{q_n(\rho)}$ reduces to a single sum, starting at $m=n$. Shifting the summation index gives the result.
\end{proof}

\begin{prop} \label{prop:G inv}
For a continuous function $g \in L^2(\R, w(x)dx)$, we have
$\mathcal F (\mathcal G g)=g$.
\end{prop}
\begin{proof}
Let $q_n$ be as in Lemma \ref{MP trans}. We show that $(\mathcal F \mathbf{q_n})(x) = p_n^{(k_2)}(x;\phi)$. From this the proposition follows, since the polynomials are dense in $L^2(\R, w(x)dx)$.

Transforming the $_2F_1$-series of $q_n$ in Lemma \ref{MP trans} by \cite[(2.3.12)]{AAR} and using the asymptotic behaviour \cite[2.3.2(16)]{Erd}, we find
\[
q_n(\rho) = \mathcal O \left(\rho^{k_1-k_2-\hf+n}\right), \qquad \rho \rightarrow \infty.
\]
From Stirling's formula it follows that $|W_0(\rho)/W_1(\rho)|= \mathcal O(e^{-\pi \rho})$ for $\rho \rightarrow \infty$. So by \eqref{assym W} and \eqref{def F} we see that $\mathcal F \mathbf{q_n}$ exists. 

We calculate
\[
I(x) = \int_0^\infty \psi_\rho(x) \overline{q_n(\rho)}W_2(\rho) d\rho,
\]
then according to Remark \ref{rem:psi} we have $\mathcal F \mathbf{q_n} = I + \overline{I}$. Writing $\psi_\rho(x)$ and $\overline{q_n(\rho)}$ as a sum, and interchanging summation and integration, gives
\[
I(x) = C \sum_{l=0}^\infty \sum_{m=0}^\infty \frac{ (1-e^{-2i\phi})^{-l} }{ l!\, (n+k_2-k_1-it+1)_l } \frac{ (k_1-ix-it)_m }{m!\, (k_1-k_2-it+1)_m (k_1+k_2-it)_m } I_{l,m},
\]
where
\[
\begin{split}
I_{l,m} =& \frac{1}{2\pi} \int_0^\infty \Ga(\hf-it+i\rho+m) \Ga(\hf-it-i\rho+m) \\ &\qquad \times  \left| \frac{ \Ga(k_1-k_2+\hf+i\rho) \Ga(n+k_2-k_1+\hf+i\rho+l) \Ga(k_1+k_2-\hf+i\rho) }{ \Ga(2i\rho) } \right|^2 d\rho,\\
C=&\, e^{-x(2\phi-\pi)}  \left( \frac{ -e^{-i\phi} }{ 1-e^{-2i\phi}} \right)^{2k_2+n}  \frac{(-1)^n}{n!}\\
&\times \frac{ \Ga(k_2-k_1+it+1)}{ \Ga(k_1+ix+it) \Ga(k_2-ix) \Ga(k_1+k_2+it) \Ga(k_2-k_1-it+n+1) \Ga(k_1-k_2-it+1)}.
\end{split}
\]
The integral $I_{l,m}$ can be evaluated by \cite[Thm.3.6.2]{AAR};
\[
\begin{split}
I_{l,m} =& \frac{ \Ga(k_2-k_1-it+n+m+l+1) \Ga( k_1-k_2-it+m+1) \Ga( k_1+k_2-it+m) }{ \Ga(k_1+k_2-it+n+m+l+1) }\\
&\times \Ga(n+l+1) \Ga(2k_2+l+n) \Ga(2k_1).
\end{split}
\]
Now the sum over $m$ is a $_2F_1$-series which can be summed by Gauss's theorem;
\[
\F{2}{1}{ k_1-ix-it, k_2-k_1-it+n+l+1}{k_1+k_2-it+n+l+1}{1} = \frac{ \Ga(k_1+ix+it) \Ga( k_1+k_2-it+n+l+1) }{ \Ga(k_2+ix+l+n+1) \Ga(2k_1) }.
\]
Then the sum over $l$ becomes a $_2F_1$-series, and after applying Euler's transformation we obtain
\[
\begin{split}
I(x) =& e^{-x(2\phi-\pi)}\frac{(-e^{i\phi})^{n-2ix}}{( 1-e^{-2i\phi} )^{ix+k_2}}   \frac{(-1)^{n}\, \Ga(2k_2+n)  }{ \Ga(k_2+ix+n+1)  \Ga(k_2-ix)} \\
&\times\F{2}{1}{k_2+ix, 1-k_2+ix}{k_2+ix+n+1}{ \frac{1}{1-e^{-2i\phi}}} .
\end{split}
\]
As in \eqref{eq:MPpol} we find from this
\[
(\mathcal F \mathbf{q_n})(x) = I(x) + \overline{I(x)} =  \frac{ (2k_2)_n}{n!}\, e^{in\phi} \F{2}{1}{-n, k_2+ix}{2k_2}{1-e^{-2i\phi}} = p_n^{(k_2)}(x;\phi).
\]
Note that the condition $\frac{\pi}{6} < \phi < \frac{5\pi}{6}$ is needed for absolute convergence of the $_2F_1$-series of $q_n$. This condition can be removed by analytic continuation.
\end{proof}

We define the Hilbert space $\mathcal M$ by
\[
\mathcal M = \overline{\text{span} \Big\{ \mathbf{q_n}\ | n \in \Z_{\geq 0} \Big\}},
\]
then $\mathcal M$ consists of functions of the form
\[
\mathbf{ f(\rho)} = \vect{f_1(\rho)}{f_2(\rho)} \in \C^2.
\]
The inner product on $\mathcal M$ is given by
\begin{equation} \label{inprod H}
\inprod{\mathbf{f}}{\mathbf{g}}_{\mathcal M} = \int_0^\infty  \vect{g_1(\rho)}{g_2(\rho)}^* 
\begin{pmatrix} 
1 & - W_0(\rho)/W_1(\rho) \\ - W_0(\rho)/\overline{W_1(\rho)} &1
\end{pmatrix}
\vect{f_1(\rho)}{f_2(\rho)}\frac{W_0(\rho)}{1- \left|\frac{W_0(\rho)}{W_1(\rho)} \right|^2} d\rho.
\end{equation}

\begin{prop} \label{prop:iso}
The operator $\mathcal F: \mathcal M \rightarrow L^2(\R, w(x)dx)$
is an isometry.
\end{prop}
\begin{proof}
Let 
\[
\mathbf{f}(\rho) = \vect{f_1(\rho)}{f_2(\rho)}, \quad \mathbf{g}(\rho) = \vect{g_1(\rho)}{g_2(\rho)}.
\]
We write out $\inprod{ \mathcal F \mathbf{f}}{\mathcal F \mathbf{g}}_{L^2(\R, w(x)dx)}$ explicitly
\[
\begin{split}
\lim_{N \rightarrow \infty} \int_{-N}^N &\int_0^\infty  \frac{W_0(\rho)}{1-\left|\frac{ W_0(\rho)}{W_1(\rho)}\right|^2 } \Big( \ph_\rho(x) f_1(\rho)-\ph_\rho^*(x) f_1(\rho) \frac{ W_0(\rho) }{\overline{W_1(\rho)}}  -\ph_\rho(x)f_2(\rho) \frac{ W_0(\rho) }{W_1(\rho)}  + \ph_\rho^*(x) f_2(\rho) \Big)  d\rho \\
\times &\overline{\int_0^\infty \frac{W_0(\si)}{1-\left|\frac{ W_0(\si)}{W_1(\si)}\right|^2 } \Big( \ph_\si(x) g_1(\si)-\ph_\si^*(x) g_1(\si) \frac{ W_0(\si) }{\overline{W_1(\si)}} -\ph_\si(x)g_2(\si) \frac{ W_0(\si) }{W_1(\si)}  + \ph_\si^*(x) g_2(\si) \Big)
d\si }\\
\times &\  w(x) dx \\
=\lim_{N \rightarrow \infty} &\int_0^\infty \int_0^\infty \frac{ W_0(\rho)W_0(\si)}{\left(1-\left|\frac{ W_0(\rho)}{W_1(\rho)}\right|^2 \right) \left( 1-\left|\frac{ W_0(\si)}{W_1(\si)}\right|^2 \right)}  \\
\times\Bigg\{ & \Big[ f_1(\rho)\overline{g_1(\si)} - f_1(\rho) \overline{g_2(\si)} \frac{ W_0(\si)}{ \overline{W_1(\si)}} - f_2(\rho)\overline{g_1(\si)} \frac{ W_0(\rho)}{ W_1(\rho)} + f_2(\rho)\overline{g_2(\si)} \frac{ W_0(\rho)W_0(\si)}{ W_1(\rho) \overline{W_1(\si)}}\Big] \inprod{\ph_\rho}{\ph_\si}_{N,N} \\
+& \Big[ f_2(\rho)\overline{g_1(\si)} -f_2(\rho) \overline{g_2(\si)} \frac{ W_0(\si)}{ \overline{W_1(\si)}} - f_1(\rho)\overline{g_1(\si)} \frac{ W_0(\rho)}{ \overline{W_1(\rho)}} + f_1(\rho)\overline{g_2(\si)} \frac{ W_0(\rho)W_0(\si)}{\overline{ W_1(\rho) W_1(\si)}}\Big] \inprod{\ph_\rho^*}{\ph_\si}_{N,N}\\
+& \Big[ f_1(\rho)\overline{g_2(\si)} -f_1(\rho) \overline{g_1(\si)} \frac{ W_0(\si)}{ W_1(\si)} - f_2(\rho)\overline{g_2(\si)} \frac{ W_0(\rho)}{W_1(\rho)} + f_2(\rho)\overline{g_1(\si)} \frac{ W_0(\rho)W_0(\si)}{ W_1(\rho) W_1(\si)}\Big] \inprod{\ph_\rho}{\ph_\si^*}_{N,N}\\
+&\Big[ f_2(\rho)\overline{g_2(\si)} - f_2(\rho)\overline{g_1(\si)} \frac{ W_0(\si)}{ W_1(\si)} - f_1(\rho)\overline{g_2(\si)} \frac{ W_0(\rho)}{ \overline{W_1(\rho)}} + f_1(\rho)\overline{g_1(\si)} \frac{ W_0(\rho)W_0(\si)}{ \overline{W_1(\rho)} W_1(\si)}\Big] \inprod{\ph_\rho^*}{\ph_\si^*}_{N,N}\Bigg\}.
\end{split}
\]
Then from Propositions \ref{prop:3}, \ref{prop:3a} and \eqref{inprod H} we obtain 
\[
\inprod{ \mathcal F \mathbf{f}}{\mathcal F \mathbf{g}}_{L^2(\R, w(x)dx)}
= \inprod{\mathbf{f}}{\mathbf{g}}_{\mathcal M}
\]
by a straightforward calculation.
\end{proof}
So far we only considered the integral transform $\mathcal F$ in the case that $\rho^2+\frac{1}{4}$ is in the continuous spectrum of the difference operator $\La$. In the next subsection we consider the discrete spectrum of $\La$.

\subsection{Discrete spectrum} \label{ssect:discrete}
From \eqref{Phi w} it follows that for $\Im(\rho)<0$ we have $\Phi_\rho(x) \in L^2(\R, w(x)dx)$. So if $c(-\rho)=0$ and $\Im(\rho)<0$, we find from Proposition \ref{prop: c-expan} that $\ph_\rho(x) = c(\rho) \Phi_\rho (x)$, and therefore $\ph_\rho(x) \in L^2(\R, w(x)dx)$. 

There are two possible cases for $c(-\rho)=0$, cf.~ Theorem \ref{thm:decomp}:
\begin{enumerate}
\item $k_2-k_1+\hf<0$, then $\rho = i(k_2-k_1+\hf+n)$, $n=0,\ldots, n_0$, where $n_0$ is the largest nonnegative
integer such that $k_2-k_1+\hf+n_0<0$,
\item $k_1+k_2-\hf<0$, then $\rho = i(k_1+k_2-\hf)$.
\end{enumerate}
Case (ii) does not occur for $k_1>1$, which is needed for convergence of the $_3F_2$-series of $\ph_\rho(x+i)$. However for $k_1 \leq 1$ we use expression \eqref{eq:Phi1} for $\Phi_\rho(x)$ (which still converges if $k_1 \leq 1$) and $\ph_\rho(x) = c(\rho) \Phi_\rho(x)$. We see that the $_3F_2$-series becomes a $_2F_1$-series of unit argument, and then, with Gauss's summation formula, we find that in case (ii) we have $\Phi_\rho(x) = e^{-2\phi x} p(x)$. Observe that case (i) and case (ii) exclude each other.

First we consider case (i). For $\rho_n = i(k_2-k_1+\hf+n)$, $0 \leq n \leq n_0$, we denote $\ph_\rho(x)$ by $\ph_{\rho_n}(x)$. We show that $\ph_{\rho_n}(x)$ is orthogonal to $\ph_{\rho_m}(x)$ and $\ph_{\rho_m}^*(x)$ for $n \neq m$. Note that $\ph_{\rho_n}$ is given by a terminating series, cf.~ \eqref{def:ph}.
\begin{prop} \label{prop:discr1}
For $m,n=0, \ldots, n_0$
\[
\begin{split}
&\inprod{ \ph_{\rho_n}}{\ph_{\rho_m}}=  \de_{nm}\left(2\pi i\, \Res{\rho = \rho_n} W_0(\rho) \right)^{-1}\\ 
&=\de_{nm}\ \frac{ e^{2t(\phi-\frac{\pi}{2})}\,\Ga(2k_2) \Ga(2k_1-2k_2-1)}{ \Ga(k_1-k_2+it) \Ga(k_1-k_2-it) \Ga(2k_1-1) } \frac{ n!\, (2k_2-2k_1+n+1)_n (2-2k_1)_n}{ (2k_2)_n (2k_2-2k_1+2)_{2n}} .
\end{split}
\]
\end{prop}
\begin{proof}
Writing out the explicit expressions \eqref{def:ph} for $\ph_{\rho_n}(x)$ and $\ph_{\rho_m}(x)$ gives
\[
\begin{split}
\int_\R \ph_{\rho_n}(x) \ph_{\rho_m}^*(x) &w(x) dx \\
=& e^{2t(\phi-\frac{\pi}{2})}\frac{1}{2\pi}\int_{-\infty}^\infty \F{3}{2}{-n, 2k_2-2k_1+n+1, k_2-ix}{2k_2, k_2-k_1+it+1}{1}\\
&\qquad \times \F{3}{2}{-m, 2k_2-2k_1+m+1, k_2+ix}{2k_2, k_2-k_1-it+1}{1} \left| \frac{ \Ga(k_2+ix)} { \Ga(k_1+it+ix) } \right|^2 dx \\
=& e^{2t(\phi-\frac{\pi}{2})}\,\sum_{k = 0}^n \sum_{l = 0}^m \frac{ (-n)_k (2k_2-2k_1+n+1)_k (-m)_l (2k_2-2k_1+m+1)_l }{ k! (2k_2)_k (k_2-k_1+it+1)_k \ l!(2k_2)_l (k_2-k_1-it+1)_l } \\ & \qquad \times\ \frac{1}{2\pi}\int_{-\infty}^\infty \frac{ \Ga(k_2-ix+k) \Ga(k_2+ix+l) }{ \Ga(k_1-it-ix) \Ga(k_1+it+ix) } dx.
\end{split}
\]
The integral inside the sum can be evaluated by \cite[\S3.3.4]{PK}
\[
\frac{1}{2\pi i} \int_{-i\infty}^{i\infty} \frac{ \Ga(a+s) \Ga(c-s) }{ \Ga(b+s) \Ga(d-s)} ds = \frac{ \Ga(a+c) \Ga(b+d-a-c-1)} { \Ga(b-a) \Ga(d-c) \Ga(b+d-1)},
\]
where $\Re(a+c-b-d)<1$ and the path of integration separates the poles of $\Ga(a+s)$ from the poles of $\Ga(c-s)$. Note that the convergence condition $2k_2-2k_1+k+l<1$ is satisfied in case (i) and $k,l \leq n_0$. Now we find for the double sum for $n \leq m$
\begin{multline*}
e^{2t(\phi-\frac{\pi}{2})}\,\frac{\Ga(2k_2) \Ga(2k_1-2k_2-1)}{ \Ga(k_1-k_2+it) \Ga(k_1-k_2-it) \Ga(2k_1-1) }\\
\times \sum_{k=0}^n \frac{ (-n)_k (2k_2-2k_1+n+1)_k }{ k! (2k_2-2k_1+2)_k } \sum_{l=0}^m \frac{ (-m)_l (2k_2-2k_1+m+1)_l (2k_2+k)_l } { l! (2k_2)_l(2k_2-2k_1+k+2)_l }.
\end{multline*}
The sum over $l$ is a terminating $_3F_2$-series, which can be evaluated by the Pfaff-Saalsch\"utz theorem
\[
\F{3}{2}{-m, 2k_2-2k_1+m+1, 2k_2+k}{2k_2, 2k_2-2k_1+k+2}{1} =
\begin{cases}
0, & k<m,\\
\displaystyle  \frac{(-1)^m m! (2-2k_1)_m }{ (2k_2)_m (2k_2-2k_1+m+2)_m }, & k=m.
\end{cases}
\]
So we find for $n \leq m$
\begin{multline*}
\int_\R \ph_{\rho_n}(x) \ph_{\rho_m}^*(x) w(x) dx = \\
\de_{nm}\ \frac{ e^{2t(\phi-\frac{\pi}{2})}\,\Ga(2k_2) \Ga(2k_1-2k_2-1)}{ \Ga(k_1-k_2+it) \Ga(k_1-k_2-it) \Ga(2k_1-1) } \frac{ n!\, (2k_2-2k_1+n+1)_n (2-2k_1)_n}{ (2k_2)_n (2k_2-2k_1+n+2)_n (2k_2-2k_1+2)_n}.
\end{multline*}
Note that from the condition $k_2-k_1+\hf+n<0$ follows that this expression is positive in case $n=m$. For $n \geq m$ we find the same result by interchanging the summations over $k$ and $l$. A straightforward calculation shows that the expression found is equal to
\[
\de_{nm}\left(2\pi i\,\Res{\rho = \rho_n} W_0(\rho) \right)^{-1}.
\]
\end{proof}

\begin{prop} \label{prop:discr1*}
For $m,n=0, \ldots, n_0$
\[
\begin{split}
\inprod{ \ph_{\rho_n}^*}{\ph_{\rho_m}}=&  \de_{nm}\left(2\pi i\, \Res{\rho = \rho_n} W_1(\rho) \right)^{-1}\\ 
=&\de_{nm}\ \frac{ e^{2t(\phi-\frac{\pi}{2})}\,\Ga(2k_2) \Ga(2k_1-2k_2-1)}{ \Ga(k_1-k_2+it) \Ga(k_1-k_2-it) \Ga(2k_1-1) }\\
& \times  \frac{(-1)^n n!(k_2-k_1+it+1)_n (2k_2-2k_1+n+1)_n (2-2k_1)_n}{(k_2-k_1-it+1)_n (2k_2)_n (2k_2-2k_1+2)_{2n}} .
\end{split}
\]
\end{prop}
\begin{proof}
From the first formula on page $142$ in \cite{AAR} we find
\[
\ph_{\rho_n}(x) = (-1)^n \frac{ (k_2-k_1-it+1)_n}{ (k_2-k_1+it+1)_n }\, \ph_{\rho_n}^*(x).
\]
Then the explicit expression follows from Proposition \ref{prop:discr1}. A straightforward calculation shows that the explicit expression is equal to the residue at $\rho=\rho_n$ of $W_1(\rho)$.
\end{proof}
A similar calculation is used for case (ii). Recall from the beginning of this subsection that in this case $\ph_{\rho_c}(x) = e^{-2\phi x} c(\rho_c) p(x)$. 
\begin{prop} \label{prop:(ii)}
Let $\rho_c=i(k_1+k_2-\hf)$, then
\[
\begin{split}
\inprod{ \ph_{\rho_c}}{\ph_{\rho_c}}&=   \left(2\pi i\,\Res{\rho = \rho_c} W_0(\rho) \right)^{-1}\\
&=e^{2t(\phi-\frac{\pi}{2})} \frac{\Ga(2k_2) \Ga(1-2k_1-2k_2)}{\Ga(1-2k_1)} \left| \frac{\Ga(k_2-k_1+it+1)}{ \Ga(k_1+k_2+it) \Ga(1-k_1-k_2+it)} \right|^2.
\end{split}
\]
\end{prop}
\begin{proof}
The proof is similar to the proof of Proposition \ref{prop:discr1}.
\end{proof}
\begin{prop} \label{prop:(ii)*}
Let $\rho_c=i(k_1+k_2-\hf)$, then
\[
\begin{split}
\inprod{ \ph_{\rho_c}^*}{\ph_{\rho_c}}&=   \left(2\pi i\,\Res{\rho = \rho_c} W_1(\rho) \right)^{-1}\\
&=e^{2t(\phi-\frac{\pi}{2})} \frac{\Ga(2k_2) \Ga(1-2k_1-2k_2)\Ga(k_2-k_1-it+1)}{\Ga(1-2k_1) \Ga(k_1+k_2-it) \Ga(1-k_1-k_2-it) \Ga(k_1-k_2+it))}.
\end{split}
\]
\end{prop}
\begin{proof}
We use Euler's reflection formula to write $p(x)$ in terms of $\Ga$-functions, then we have 
\[
\inprod{ \ph_{\rho_c}^*}{\ph_{\rho_c}}= e^{2t(\phi-\frac{\pi}{2})} \overline{c(\rho_c)}^2 \frac{1}{2\pi}\int_\R \frac{ \Ga(k_2+ix) \Ga(k_2-ix) \Ga(k_1-it-ix) }{ \Ga(k_1+it+ix) \Ga(1-k_1-it-ix) \Ga(1-k_1-it-ix)} dx.
\]
We use a special case of \cite[(4.5.1.2)]{Sl1} to evaluate the integral;
\[
\begin{split}
\frac{1}{2\pi i} &\int_{-i\infty}^{i\infty} \frac{ \Ga(a+s) \Ga(b-s) \Ga(c-s) }{ \Ga(d+s) \Ga(e-s) \Ga(f-s) } ds =\\
 &\frac{ \Ga(a+b) \Ga(a+c) }{ \Ga(d-a) \Ga(a+e) \Ga(a+f) } \F{3}{2}{ a+b, a+c, 1+a-d}{ a+e, a+f }{1},
\end{split}
\]
where $\Re(d+e+f-a-b-c)>0$ and the path of integration separates the poles of $\Ga(a+s)$ from the poles of $\Ga(b-s)$ and $\Ga(c-s)$. We put 
\begin{gather*}
s= -k_2+ix, \quad a= 2k_2, \quad  b = 0, \quad  c = k_1-k_2-it, \\
d = k_1+k_2+it, \quad e = 1-k_1-k_2-it, \quad f= 1-k_1-k_2-it,
\end{gather*}
then we find that the $_3F_2$-series reduces to a $_2F_1$-series, which can be evaluated by Gauss's summation formula. From this the result follows.
\end{proof}

Next we show that if $\rho^2+\frac{1}{4}$ is in the discrete spectrum of $\La$, $\ph_\rho(x)$ is orthogonal to $\ph_\si(x)$ and $\ph_\si^*(x)$, if $\si$ is real. 
\begin{prop} \label{prop:orth}
For $\si \in [0,\infty)$ and $\rho= \rho_n$, or $\rho=\rho_c$, we have
\[
\inprod{ \ph_\rho}{\ph_\si} =0, \qquad \inprod{ \ph_\rho}{\ph_\si^*} =0.
\]
\end{prop}
\begin{proof}
By Propositions \ref{prop: int=wrons} and \ref{prop:eigenf}
\[
\lim_{N,M \rightarrow \infty} \int_{-M}^N \ph_\rho(x) \ph_\si^*(x) w(x) dx =\lim_{N,M \rightarrow \infty} \frac{[\ph_\rho,\ph_\si](N) - [\ph_\rho, \ph_\si](-M)}{\rho^2-\si^2}.
\]
We show that the limit of each Wronskian is zero. Let
\[
f(x)=\Big(\ph_\rho(x) \ph^*_\si(x-i) - \ph_\rho(x-i) \ph^*_\si(x) \Big) \al_-(x) w(x).
\]
First we consider the asymptotic behaviour of $f(-M)$ and $f(N)$.

For $\rho=\rho_n$ and $0 \leq y \leq 1$ we find from Lemmas \ref{eq:asymp w}, \ref{lem: assympt} and Proposition \ref{prop: c-expan} 
\[
f(x+iy) = \mathcal O(|x|^{\hf-k_1+k_2-n}), \qquad x \rightarrow \pm \infty.
\]
For $\rho=\rho_c$ we find
\[
f(x+iy) = \mathcal O(|x|^{2k_1+2k_2-1}), \qquad x \rightarrow \pm \infty.
\]
Here the implied constants do not depend on $y$. Then dominated convergence gives the result.

The proof for $\inprod{ \ph_\rho}{\ph_\si^*} =0$ runs along the same lines.
\end{proof}

\begin{rem}
The explicit calculations in this subsection can be carried out because of the choice \eqref{def: periodic} of the $i$-periodic function $p(x)$. It is not likely that with another choice for the function $p(x)$ all the calculations can be done explicitly. This gives another (heuristic) argument for the choice of $p(x)$.
\end{rem}

\subsection{The continuous Hahn integral transform} We combine the results of subsection \ref{ssect: continuous} with the results of subsection \ref{ssect:discrete}.

Let $k_1, k_2>0$, $0<\phi<\pi$ and $t \in \R$. Let $\ph_\rho(x)$ be the function given by
\[
\ph_\rho(x) =\frac{ 
e^{-x(2\phi-\pi)}}{\Ga(k_1+ix+it)\Ga(k_1-ix-it)} \F{3}{2}{ k_2-ix,
k_2-k_1+\hf+i\rho, k_2-k_1+\hf-i\rho} {2k_2, k_2-k_1+it+1} {1}.
\]
We denote
\[
\begin{split}
h(\rho) &= \frac{W_0(\rho)}{W_1(\rho)} = \frac{\Ga(\hf+it+i\rho) \Ga( \hf + it -i\rho) }{ \Ga(k_1-k_2+it) \Ga(k_2-k_1+it+1) },\\
W(\rho) &= \frac{W_0(\rho)}{1-|h(\rho)|^2} = \frac{ 1}{2\pi} e^{- t(2\phi-\pi)}\, \left| \frac{ \Ga( k_2-k_1+\hf+i\rho)^2 \Ga( k_1-k_2+\hf+i\rho) \Ga(k_1+k_2-\hf+i\rho) }{ \Ga(2k_2) \Ga(k_2-k_1+it+1) \Ga(2i\rho) }\right|^2,\\
w(x) &=\frac{1}{2\pi} e^{(2x+t)(2\phi-\pi)} | \Ga(k_1+it+ix) \Ga(k_2+ix) |^2.
\end{split}
\]
Let $\mathcal M$ be the Hilbert space given by
\[
\mathcal M = \overline{\text{span} \Big\{ \mathbf{q_m}\ | \ m \in \Z_{\geq 0} \Big\} },
\] 
where $\mathbf{q_m}$ is as in Lemma \ref{MP trans}. 
For functions $\mathbf{f} = \vect{f_1}{f_2}$, the inner product on $\mathcal M$ is given by:\\

(i) For $k_1+k_1-\hf \geq 0$, $k_2-k_1+\hf \geq 0$
\[
\inprod{\mathbf{f}}{\mathbf{g}}_{\mathcal M}= \int_0^\infty \vect{g_1(\rho)}{g_2(\rho)}^* \begin{pmatrix}
1 & -h(\rho) \\
-\overline{h(\rho)} &1
\end{pmatrix}
\vect{f_1(\rho)}{f_2(\rho)}W(\rho)d\rho. 
\]

(ii) For $k_1+k_2-\hf<0$, $\rho_c=i(k_1+k_2-\hf)$, 
\[
\begin{split}
\inprod{\mathbf{f}}{\mathbf{g}}_{\mathcal M}= \int_0^\infty \vect{g_1(\rho)}{g_2(\rho)}^* &\begin{pmatrix}
1 & -h(\rho) \\
-\overline{h(\rho)} &1
\end{pmatrix}
\vect{f_1(\rho)}{f_2(\rho)} W(\rho)d\rho\\ & + \vect{g_1(\rho_c)}{g_2(\rho_c)}^* \begin{pmatrix}
 1& -h(\rho_c) \\
-\overline{h(\rho_c)} & 1
\end{pmatrix}
\vect{f_1(\rho_c)}{f_2(\rho_c)}\, 2\pi i\, \Res{\rho=\rho_c}W(\rho). 
\end{split}
\]

(iii) For $k_2-k_1+\hf<0$, $\rho_n = i(k_2-k_1+\hf+n)$, $n=0, \ldots, n_0$, where $n_0$ is the largest integer such that $-i\rho_{n_0}<0$,
\[
\begin{split}
\inprod{\mathbf{f}}{\mathbf{g}}_{\mathcal M}= \int_0^\infty & \vect{g_1(\rho)}{g_2(\rho)}^* \begin{pmatrix}
1 & -h(\rho) \\
-\overline{h(\rho)} & 1
\end{pmatrix}
\vect{f_1(\rho)}{f_2(\rho)}W(\rho)d\rho\\
&+\sum_{n=0}^{n_0} \vect{g_1(\rho_n)}{g_2(\rho_n)} ^*
\vect{f_1(\rho_n)}{f_2(\rho_n)}\, \pi i\, \Res{\rho=\rho_n} W(\rho)(1-|h(\rho)|^2). 
\end{split}
\]
Observe that for $m > n$ we have $\mathbf{q_m}(\rho_n)=\mathbf{0}$. For $m \leq n$ it follows from the way $q_m$ is calculated in Proposition \ref{MP trans}, that
$q_m(\rho_n) = \overline{h(\rho_n)} \overline{q_m(\rho_n)}$. 

For a continuous function $\mathbf{f}\in \mathcal M$ we define the linear operator $\mathcal F: \mathcal M \rightarrow L^2(\R, w(x)dx)$ by 
\[
(\mathcal F \mathbf{f})(x) = \left\langle \mathbf{f} ,\vect{\ph_\rho^*(x)}{\ph_\rho(x)} \right\rangle_{\mathcal M}.
\]
We call $\mathcal F$ the continuous Hahn integral transform.
\begin{thm} \label{thm:CH-transform}
The continuous Hahn integral transfrom $\mathcal F:\mathcal M \rightarrow L^2(\R, w(x)dx)$ is unitary and its inverse is given by
\[
\left(\mathcal F^{-1} g \right)(\rho)=\int_\R g(x) \vect{\ph_\rho^*(x)}{\ph_\rho(x)} w(x)dx .
\]
\end{thm}
\begin{proof}
For case (i) this follows from Propositions \ref{prop:F inv}, \ref{prop:G inv} and \ref{prop:iso}.

For case (ii) we only have to check that Propositions \ref{prop:3} and \ref{prop:3a} still hold with the discrete mass point in $\rho=\rho_c$ added to the integral. From Propositions \ref{prop:(ii)}, \ref{prop:(ii)*} and \ref{prop:orth} we find
\[
\begin{split}
f(\rho_c) \inprod{\ph_{\rho_c}}{ \ph_\si}=& 
\begin{cases}
0, & \si \in [0,\infty),\\
\left(2\pi i \Res{\si=\rho_c} W_0(\si) \right)^{-1}f(\rho_c), & \si=\rho_c, 
\end{cases}\\
f(\rho_c) \inprod{\ph_{\rho_c}^*}{ \ph_\si}= &
\begin{cases}
0, & \si \in [0,\infty),\\
\left(2\pi i \Res{\si=\rho_c} W_1(\si) \right)^{-1}f(\rho_c), & \si=\rho_c. 
\end{cases}
\end{split}
\]
Now the proof for case (ii) is completely analogous to the proof of case (i).

For case (iii) injectivity and surjectivity of $\mathcal F$ can be proved in the same way as case (i). We check that $\mathcal F$ is an isometry. The continuous part follows from Proposition \ref{prop:iso}, so we only have to check for the discrete part. We write out $\inprod{\mathcal F \mathbf{q_k}}{\mathcal F \mathbf{q_l}}_{L^2(\R, w(x)dx}$, $k,l \in \Z_{\geq 0}$, for the discrete part of $\mathcal F \mathbf{q_k}$ and $\mathcal F \mathbf{q_l}$. From Propositions \ref{prop:discr1}, \ref{prop:discr1*} and $q_k(\rho_n) = \overline{h(\rho_n)} \overline{q_k(\rho_n)}$ we find, for $n,m =0, \ldots, n_0$,
\[
\begin{split}
\int_\R &\ \sum_{n=0}^{n_0}\Big( \ph_{\rho_n}(x) \overline{q_k(\rho_n)} + \ph_{\rho_n}^*(x) q_k(\rho_n) \Big) \, \pi i\, \Res{\rho=\rho_n} W_0(\rho)\\
&\times \sum_{m=0}^{n_0}\Big( \ph_{\rho_m}(x) \overline{q_l(\rho_m)}  + \ph_{\rho_m}^*(x) q_l(\rho_m)\Big)\, \pi i\, \Res{\rho=\rho_m} W_0(\rho)\ w(x) dx\\
=&\sum_{n,m=0}^{n_0} \Big( \overline{ q_k(\rho_n) q_l(\rho_m)}   \inprod{\ph_{\rho_n}}{\ph_{\rho_m}^*} +
\overline{ q_k(\rho_n)} q_l(\rho_m) \inprod{\ph_{\rho_n}}{\ph_{\rho_m}} 
+ q_k(\rho_n) \overline{q_l(\rho_m)} \inprod{\ph_{\rho_n}^*} {\ph_{\rho_m}^*}\\ 
&\qquad+  q_k(\rho_n) q_l(\rho_m) \inprod{\ph_{\rho_n}^*}{\ph_{\rho_m}} \Big)\, (\pi i)^2\, \Res{\rho=\rho_n} W_0(\rho) \,  \Res{\rho=\rho_m} W_0(\rho)\\
= &\sum_{n=0}^{n_0}\Big( q_k(\rho_n) \overline{q_l(\rho_n)}  + \overline{ q_k(\rho_n)} q_l(\rho_n) \Big)\, \pi i\, \Res{\rho=\rho_n} W_0(\rho). 
\end{split}
\]
Here we recognize the discrete part of the inner product $\inprod{\mathbf{q_k}}{\mathbf{q_l}}_{\mathcal M}$. Combined with Propositions \ref{prop:iso} and \ref{prop:orth} this shows that $\mathcal F$ acts isometric on the basis elements $\mathbf{q_k}$. By linearity $\mathcal F$ extends to an isometry. 
\end{proof}

The continuous Hahn integral transform in case (i) corresponds exactly to the integral transform we found in \S\ref{ssec:formal} by formal computations.\\

\begin{rem} \label{rem: diff eig}
Let us denote the operator $\La$ by $\La(k_1,k_2,t)$, let $w(x) = w(x;k_1,k_2,t)$, and let $T_t$ denote the shift operator. Observe that $w(x+t;k_2,k_1,-t)=w(x;k_1,k_2,t)$. It is clear that $L^2(\R, w(x;k_1,k_2,t)dx)$ is invariant under the action of $T_t$. A short calculation shows that $T_{-t} \circ \La(k_1,k_2,t) \circ T_t = \La(k_2,k_1,-t)$, so  $\ph_\rho(x+t;-t,k_2,k_1,\phi)$ is an eigenfunction of $\La(k_1,k_2,t)$ for eigenvalue $\rho^2+\frac{1}{4}$. Going through the whole machinery of this section again, then gives another spectral measure of $\La(k_1,k_2,t)$, namely the one we found with $k_1 \mapsto k_2$, $k_2 \mapsto k_1$ and $t \mapsto -t$. 
\end{rem}

Finally we compare the spectrum of the difference operator $\La$ with the tensor product decomposition in Theorem \ref{thm:decomp}. The discrete term in case (ii) in this section corresponds to one complementary series representation in the tensor product decomposition in Theorem \ref{thm:decomp}. Case (iii) does not occur in Theorem \ref{thm:decomp}, since it is assumed that $k_1 \leq k_2$. If we no longer assume this, the discrete terms in case (iii) correspond to a finite number of positive discrete series representations. 

The finite number of negative discrete series in Theorem \ref{thm:decomp} can be obtained as described in Remark \ref{rem: diff eig}. In order to obtain the Clebsch-Gordan coefficients in this case from the summation formula in Theorem \ref{Thm:sum}, we need to consider different overlap coefficients for the continuous series representations. Let $v_n(x;\la, \eps,\phi) := u_{-n}(x;\la, -\eps,\pi -\phi)$, where $u_n$ is the Meixner-Pollaczek function as defined by \eqref{def:MPfunctions}. Then 
\[
\vect{\tilde v_{\rho, \eps}(x)}{\tilde v_{\rho, \eps}^*(x)} = \sum_{n=-\infty}^\infty \vect{v_n(x;-\hf+i\rho, \eps,\pi-\phi) }{v_n^*(x;-\hf+i\rho, \eps,\pi-\phi) }\,e_n,
\]
is a generalized eigenvector of $\pi^{\rho,\eps}(X_\phi)$, see \cite[\S4.4.11]{Koe}. From Theorem \ref{Thm:sum}, with $(k_1,k_2,x_1,x_2,p) \mapsto (k_2, k_1, x_2,x_1,-p)$ we find the Clebsch-Gordan coefficients for the eigenvector $\vect{\tilde v_{\rho,\eps}(x_1-x_2)}{\tilde v_{\rho,\eps}^*(x_1-x_2)}$, and these Clebsch-Gordan coefficients are multiples of continuous Hahn functions $\ph_\rho(x_1;x_2-x_1,k_2,k_1,\phi)$ and $\ph_\rho^*(x_1;x_2-x_1,k_2,k_1,\phi)$. In this case the discrete mass points in the measure for the continuous Hahn transform correspond to one complementary series representation, or a finite number of negative discrete series representations

\end{document}